\documentclass[11pt]{amsart}

\theoremstyle{plain}
\newtheorem{theorem}{Theorem}[section]
\newtheorem{defi}{Definition}[section]

\newtheorem{lemma}{Lemma}[section]
\newtheorem{proposition}{Proposition}[section]
\numberwithin{equation}{section}
\newcommand{\dem}{\medskip \par \noindent \mbox{\bf Proof. }}
\def\ep{\hfill{$\Box $}}

\begin{document}

\title[Function  spaces for the Gaussian measure]{Some results on Gaussian Besov-Lipschitz spaces and Gaussian Triebel-Lizorkin spaces.}

\author{Ebner Pineda}
\address{Departamento de Matem\'{a}tica,  Decanato de Ciencia y Tecnologia, UCLA
 Apartado 400 Barquisimeto 3001 Venezuela.}
\email{[Ebner Pineda]epineda@uicm.ucla.edu.ve}
\author{Wilfredo O. Urbina R.}
\address{Departamento de Matem\'{a}ticas, Facultad de Ciencias, UCV, Los Chaguaramos, Caracas 1041 Venezuela, and
 Department of Mathematical Sciences, DePaul University  Chicago, Il, 60614, USA.}
\email{[Wilfredo Urbina]wilfredo.urbina@ciens.ucv.ve, wurbina@depaul.edu}
\thanks{\emph{2000 Mathematics Subject Clasification} Primary 42C10; Secondary 26A24}
\thanks{$^{(1)}$Partially supported by Grant FONACIT \#G-97000668.}
\maketitle

\begin{abstract}
In this paper we define  Besov-Lipschitz
 and  Triebel-Lizorkin 
 spaces in the context of Gaussian harmonic analysis, 
 the harmonic analysis of Hermite polynomial expansions.
We study inclusion relations among them, some interpolation results and continuity 
results of some important operators (the Ornstein-Uhlenbeck and 
 the Poisson-Hermite semigroups and the Bessel potentials) on them. We also prove that the Gaussian
Sobolev spaces $L^p_\alpha(\gamma_d)$ are contained in them.
The proofs are general enough to allow extensions of these results to the case of Laguerre or Jacobi expansions and even further in the general framework of diffusions semigroups.
\end{abstract}

\section{Introduction}

Let us consider the Gaussian measure $\gamma_d(x)=\frac{e^{-\left|x\right|^2}%
}{\pi^{d/2}}$ with $x\in\mathbb{R}^d$ and the Ornstein-Uhlenbeck
differential operator
\begin{equation}\label{OUop}
L=\frac12\triangle_x-\left\langle x,\nabla _x\right\rangle.
\end{equation}

Let $\beta=(\beta _1,...,\beta_d)\in\mathbb{N}^d$ be a
multi-index, let $\beta
!=\prod_{i=1}^d\beta _i!,$ $\left| \beta \right| =\sum_{i=1}^d\beta _i,$ $%
\partial _i=\frac \partial {\partial x_i},$ for each $1\leq i\leq d$ and $%
\partial ^\beta =\partial _1^{\beta _1}...\partial _d^{\beta _d}.$

Let us consider the normalized Hermite polynomials of order
$\beta$, in $d$ variables
\begin{equation}
h_\beta (x)=\frac 1{\left( 2^{\left| \beta \right| }\beta
!\right)
^{1/2}}\prod_{i=1}^d(-1)^{\beta _i}e^{x_i^2}\frac{\partial ^{\beta _i}}{%
\partial x_i^{\beta _i}}(e^{-x_i^2}),
\end{equation}
then, it is well known, that the Hermite polynomials are
eigenfunctions of $L$,
\begin{equation}\label{eigen}
L h_{\beta}(x)=-\left|\beta \right|h_\beta(x).
\end{equation}
Given a function $f$ $\in L^1(\gamma _d)$ its
$\beta$-Fourier-Hermite coefficient is defined by
\[
\hat{f}(\beta) =<f, h_\beta>_{\gamma_d}
=\int_{\mathbb{R}^d}f(x)h_\beta (x)\gamma _d(dx).
\]
Let $C_n$ be the closed subspace of $L^2(\gamma_d)$ generated by
the linear combinations of $\left\{ h_\beta \ :\left| \beta
\right| =n\right\}$. By the orthogonality of the Hermite
polynomials with respect to $\gamma_d$ it is easy to see that
$\{C_n\}$ is an orthogonal decomposition of $L^2(\gamma_d)$,
$$ L^2(\gamma_d) = \bigoplus_{n=0}^{\infty} C_n$$
which is called the Wiener chaos.

Let $J_n$ be the orthogonal projection  of $L^2(\gamma_d)$ onto
$C_n$, then if $f\in L^2(\gamma_d)$
\[
J_n f=\sum_{\left|\beta\right|=n}\hat{f}(\beta) h_\beta.
\]
Let us define the Ornstein-Uhlenbeck semigroup $\left\{
T_t\right\} _{t\geq 0}$ as
\begin{eqnarray}\label{01}
\nonumber T_t f(x)&=&\frac 1{\left( 1-e^{-2t}\right) ^{d/2}}\int_{\mathbb{R}^d}e^{-\frac{%
e^{-2t}(\left| x\right| ^2+\left| y\right| ^2)-2e^{-t}\left\langle
x,y\right\rangle }{1-e^{-2t}}}f(y)\gamma _d(dy)\\
& = & \frac{1}{\pi^{d/2}(1-e^{-2t})^{d/2}}\int_{\mathbb R^d} e^{-
\frac{|y-e^{-t}x|^2}{1-e^{-2t}}} f(y) dy
\end{eqnarray}
The family $\left\{ T_t\right\}_{t\geq 0}$ is a strongly
continuous Markov semigroup on $L^p(\gamma_d)$, $1 \leq p <
\infty$, with infinitesimal generator $L$. Also, by a change of
variable we can write,
\begin{equation}\label{t1}
T_t f(x)=\int_{\mathbb{R}^d} f(\sqrt{1-e^{-2t}}u + e^{-t}x)\gamma
_d(du).
\end{equation}

Now, by Bochner subordination formula, see Stein \cite{se70}, we
define the Poisson-Hermite semigroup $\left\{ P_t\right\} _{t\geq
0}$ as
\begin{equation}\label{02}
\nonumber P_t f(x)=\frac 1{\sqrt{\pi }}\int_0^{\infty} \frac{e^{-u}}{\sqrt{u}}T_{t^2/4u}f(x)du
=\int_0^{\infty} T_s f(x) \mu^{(1/2)}_t(ds),
\end{equation}
where the measure
\begin{equation}\label{onesided1/2}
\mu^{(1/2)}_t(ds) = \frac t{2\sqrt{\pi
}}\frac{e^{-t^2/4s}}{s^{3/2}}ds = g(t,s) ds,
\end{equation}
is called the one-side stable measure on $(0, \infty)$ of order
$1/2$.

The family $\left\{ P_t\right\}_{t\geq 0}$ is also a strongly
continuous semigroup on $L^p(\gamma_d)$, $1 \leq p < \infty$,
with infinitesimal generator $-(-L)^{1/2}$. From (\ref{01}) we
obtain, after the change of variable $r=e^{-t^2/4u}$,
\begin{eqnarray}\label{03}
\nonumber P_t f(x)&=&\frac 1{2\pi
^{(d+1)/2}}\int_{\mathbb{R}^d}\int_0^1t\frac{\exp \left( t^2/4\log
r\right) }{(-\log r)^{3/2}}\frac{\exp \left( \frac{-\left|
y-rx\right| ^2}{1-r^2}\right) }{(1-r^2)^{d/2}}\frac{dr}rf(y)dy\\
&=& \int_{\mathbb{R}^d} p(t,x,y) f(y)dy,
\end{eqnarray}
with
\begin{equation}
p(t,x,y) = \frac 1{2\pi ^{(d+1)/2}}\int_0^1t\frac{\exp \left(
t^2/4\log r\right) }{(-\log r)^{3/2}}\frac{\exp \left(
\frac{-\left| y-rx\right| ^2}{1-r^2}\right)
}{(1-r^2)^{d/2}}\frac{dr}r.
\end{equation}
In what follows, we will often going to use the notation $$u(x,t) = P_{t}f(x)$$.

Observe that by (\ref{eigen}) we have that
\begin{equation}\label{OUHerm}
T_t h_\beta(x)=e^{-t\left| \beta\right|}h_\beta(x),
\end{equation}
and
\begin{equation}\label{PHHerm}
 P_t h_\beta(x)=e^{-t\sqrt{\left| \beta\right|}}h_\beta(x)
\end{equation}

Let us observe that since $|| T_t f -f ||_{p, \gamma_d} \rightarrow 0$ and $|| P_t f -f ||_{p, \gamma_d} \rightarrow 0$ as $ t \rightarrow 0$ then $\{T_t\}$   and $\{P_t\}$  play the role of ``approximation of the identity" in Gaussian setting. Moreover they are, up to now, the only  approximations of identity known. Therefore following H. Triebel, see  \cite{trie2} section 2.6.4 Harmonic and Thermic extensions pag 152, we are going to use them to define Gaussian Besov-Lipschitz
$B_{p,q}^{\alpha}(\gamma_d)$ and Gaussian Triebel-Lizorkin  $
F_{p,q}^{\alpha}(\gamma_d)$ spaces. An open problem then  is to find alternative definitions of those spaces and give a more explicit description about the type of regularity that they actually describe.

On the other hand, the possibility of characterize the Gaussian Besov-Lipchitz spaces in terms  of  modulus of smoothness, as it is done in the classical case, would be possible only if the classical  translation operator $\tau_yf(x) = f(x+y)$ is replaced  for a more suitable  translation operator since  the spaces $L^p(\gamma_d)$ are not in general closed under the action of $\tau_y$, for instance, in the one dimensional case, let us take the function $f(x) =e^{|x|^2-|x|}$, then it is clear that $f \in L^1(\gamma_1)$
but it is easy to see that $\tau_1f(x) = f(x+1)=e^{|x+1|^2-|x+1|} \notin L^1(\gamma_1)$. This point requieres further investigations.\\

 For $\alpha>0$, the Fractional Integral or Riesz potential of order $\alpha$, $I_\alpha^{\gamma}$, with respect to the Gaussian measure is defined formally as
\begin{equation}\label{i1}
I_\alpha^{\gamma}=(-L)^{-\alpha/2}\Pi_{0},
\end{equation}
where,
$\Pi_{0}f=f-\displaystyle\int_{\mathbb{R}^{d}}f(y)\gamma_{d}(dy)$,
for $f\in L^{2}(\gamma_{d})$. That means that for  the Hermite
polynomials $\{h_\beta\}$, for $\left|\beta\right|>0$,
\begin{equation}\label{e4}
I_\alpha^{\gamma}h_\beta(x)=\frac 1{\left|
\beta\right|^{\alpha/2}}h_\beta(x),
\end{equation}
and for $\beta=\overline{0}, \,
I_{\alpha}^{\gamma}(h_{\overline{0}})=0.$ Then by linearity can be extended to any polynomial.
It is easy to see that if $f$ is a polynomial,
\begin{equation}\label{e3}
I_\alpha^{\gamma} f(x)  =\frac 1{\Gamma(\alpha)}\int_0^{\infty}
t^{\alpha-1}(P_t f(x) -P_{\infty} f(x))\,dt.
\end{equation}
Moreover  by P. A. Meyer's multiplier theorem, see \cite{me}, $I_\alpha^{\gamma}$ admits a continuous extension  to $ L^p(\gamma _d)$,  $1 < p < \infty$, and (\ref{e3}) can be extended for $f \in L^{p}(\gamma_d)$, see \cite{sp97}. Also if $f$ $\in C_B^2(\mathbb{R}^d)$ such that $\int_{\mathbb{R}^d}f(y)\gamma_d(dy)=0$, then
\begin{equation}\label{p2}
I_\alpha^{\gamma} f =-\frac {1}{\alpha \Gamma(\alpha)}\int_0^{\infty} t^\alpha \frac{\partial}{\partial t}P_t f dt,
\end{equation}
see \cite{lu}.

The Bessel Potential of order $\alpha>0,$
$\mathcal{J}_\alpha^{\gamma}$, associated to the Gaussian  measure
is defined formally as
\begin{eqnarray}
\mathcal{J}_\alpha^{\gamma}  = (I-L)^{-\alpha /2},
\end{eqnarray}
meaning that for the Hermite polynomials we have,
\begin{eqnarray*}
\mathcal{J}^\gamma_\alpha h_\beta(x)=\frac 1{(1+\left|
\beta\right|)^{\alpha/2}}h_\beta(x).
\end{eqnarray*}
Again  by linearity can be extended to any polynomial and Meyer's theorem allows us to extend Bessel Potentials
to a \mbox{continuous} operator on $L^p(\gamma_d),$ $1 < p <
\infty$. Additionaly, it is easy to see that $\mathcal{J}^\gamma_\alpha$ is a bijection over the set of polynomials ${\mathcal P}$.The Bessel potentials can be
represented as
\begin{equation}\label{Beselrepre}
\mathcal{J}^\gamma_\alpha
f(x)=\frac{1}{\Gamma(\alpha)}\int_{0}^{+\infty}t^{\alpha}e^{-t}P_{t}f(x)\frac{dt}{t},
\end{equation}
for more details see \cite{forscotur}.
Moreover $\{\mathcal{J}_\alpha^{\gamma}\}_\alpha$
is a strongly continuous semigroup on $L^p(\gamma_d)$, $1 \leq p <\infty$, with infinitesimal generator $\frac{1}{2}\log (I-L).$

The fractional derivate of order $\alpha >0$ with respect to the
Gaussian measure $D_\alpha^{\gamma}$, is defined  formally as
\begin{equation}
D_\alpha^{\gamma}=(-L)^{\alpha/2},
\end{equation}
meaning that for the Hermite polynomials, we have
\begin{equation}\label{e6}
D_\alpha^{\gamma}h_\beta(x)=\left| \beta\right|^{\alpha/2}
h_\beta(x),
\end{equation}
 thus by linearity can be extended to any polynomial.
 
  The fractional derivate $D_\alpha^{\gamma}$  with respect to the Gaussian
measure was first introduced in \cite{lu}. For more detail we
refer to that article. Also see \cite{ebner} for improved  and
simpler proofs of some results contained there.

Now, if $f$ is a polynomial, by the linearity of the operators
$I_\alpha^{\gamma}$ and $D_\alpha^{\gamma}$, (\ref{e4}) and
(\ref{e6}), we get
\begin{equation}\label{di}
 \Pi_0f=I_\alpha^{\gamma}(D_\alpha^{\gamma} f)=D_\alpha^{\gamma}(I_\alpha^{\gamma} f).
\end{equation}

The Gaussian Sobolev spaces of order $\alpha \geq 0$,
$L_{\alpha}^p(\gamma_d), 1<p<\infty,$ can be obtained, as in the classical case, as the image
of $L^p(\gamma_d)$ under the Bessel potential
$\mathcal{J}_\alpha^{\gamma}$, with the norm
\begin{equation}\label{e12}
\left\|f \right\|_{p,\alpha}:=\left\| (I-L)^{\alpha/2}f
\right\|_{p,\gamma_d}.
\end{equation}
Also they can be defined as the completion of the set of polynomials ${\mathcal P}$ with respect to that norm, see
\cite{wat} and therefore ${\mathcal P}$ is trivially dense there. Let us remember that it can be proved that the set of polynomials ${\mathcal P}$ is also dense in $L^p(\gamma_d), 1<p<\infty $, see \cite{berg}.
The fractional derivative $D_\alpha^{\gamma}$ can be used to  characterize the Gaussian Sobolev spaces $L^p_\alpha(\gamma_d)$ see \cite{lu} .

As usual in what follows $C$ represents a constant that is not necessarily the same in each occurrence.

We wish to express our thanks to Prof. A. Eduardo Gatto for his
useful conversations,  and suggestions. We also want to thanks
Prof. Hugo Aimar for an important observation that lead us to
Theorem 2.4. Also we want to thank the referees for all their suggestions, comments and observations 
which have improved not only the presentation of our paper but also some of the results  obtained.

\section{The  main results}
As it was already mentioned in the introduction, the main objective of this
paper is to introduce the Gaussian Besov-Lipschitz
$B_{p,q}^{\alpha}(\gamma_d)$ and the Gaussian Triebel-Lizorkin  $
F_{p,q}^{\alpha}(\gamma_d)$ spaces, for any $\alpha \geq 0$. We will follow E. Stein \cite{se70} scheme
to define and study the $B_{p,q}^{\alpha}(\gamma_d)$ spaces, but since the Poisson-Hermite semigroup is not a convolution semigroup the proofs of the results will be totally different to the ones in Stein's book. We will use, in an essential way, the representation of the Poisson-Hermite semigroup  (\ref{02}) using the one-side stable measure, $\mu^{(1/2)}_t$ defined in (\ref{onesided1/2}). From that fact, it is then clear that similar constructions are possible for the harmonic analysis of Laguerre or Jacobi polynomial expansions and even further in the framework of general diffusion semigroups but we are not going to consider those cases here. Let us point out that Hermite, Laguerre and Jacobi are the only cases of diffusion semigroups associated to orthogonal polynomials, see Mazet \cite{ma}.

On the other hand, Besov-Lipschitz spaces can be also obtained as interpolated spaces using interpolation theory for semigroups defined on a Banach space, see for instance Chapter 3 of \cite{ButzBer} or \cite{trie}.

We will need some technical results for the measure $\mu^{(1/2)}_t$. First, in what follows since $\mu^{(1/2)}_t(ds) = \frac t{2\sqrt{\pi}}\frac{e^{-t^2/4s}}{s^{3/2}}ds = g(t,s) ds$,  for any  $k\in \mathbb{N}$,
the notation $\frac{\partial^{k}}{\partial t^{k}}\mu_{t}^{(1/2)}(ds)$ will denote 
 \begin{equation}
\frac{\partial^{k}}{\partial t^{k}}\mu_{t}^{(1/2)}(ds) := \frac{\partial^{k} g(t,s)}{\partial t^{k}}ds.
\end{equation}
Then by induction it can be seen that
\begin{equation}\label{inductdient}
 \frac{\partial^{k}\mu_{t}^{(1/2)}}{\partial t^{k}}(ds)=\big(\sum_{\begin{array}{c}i\in
\mathbb{Z}, j\in \mathbb{N},\\  0\leq j\leq k, 2j-i=k\end{array}}
a_{i,j}\frac{t^{i}}{s^{j}}\big)\mu_{t}^{(1/2)}(ds)
\end{equation}
where $\{a_{i,j}\}$ is a (finite) set of  constants.

Moreover, using the change of variable $\displaystyle
u=\frac{t^{2}}{4s}$, it is easy to see that given $k\in \mathbb{N}$ and $t>0$
\begin{equation}\label{intoneside}
\int_{0}^{+\infty}\frac{1}{s^{k}}\mu_{t}^{\frac{1}{2}}(ds)=\frac{C_{k}}{t^{2k}},
\end{equation}
where
$C_{k}=\frac{2^{2k}\Gamma(k+\frac{1}{2})}{\pi^{\frac{1}{2}}}.$
Finally, using the two previous results we get that if 
 $k\in \mathbb{N}$ and $t>0$, then
\begin{equation}\label{onesideineq}
 \int_{0}^{+\infty}| \frac{\partial^{k}}{\partial
t^{k}}\mu_{t}^{(1/2)}|(ds)\leq\frac{C_{k}}{t^{k}}.
\end{equation}

Now, considering the maximal function of the
Ornstein-Uhlenbeck semigroup,
$$T^{\ast}f(x)=\displaystyle\sup_{t>0}|T_{t}f(x)|,$$
 we have the following inequality that will be used later,

\begin{lemma}\label{lematec5}
$$
|\frac{\partial^{k}P_{t}f(x)}{\partial t^{k}}| \leq C_{k} \,
T^{\ast}f(x) t^{-k}.
$$
\end{lemma}
\dem
Using (\ref{onesideineq}) and the dominated convergence theorem, we have
\begin{eqnarray*}
|\frac{\partial^{k}P_{t}f(x)}{\partial
t^{k}}|&=&|\int_{0}^{+\infty}T_{s}f(x) \frac{\partial^{k}
}{\partial t^{k}}\mu_{t}^{(1/2)}(ds)|
\leq\int_{0}^{+\infty}|T_{s}f(x)| | \frac{\partial^{k}
}{\partial
t^{k}}\mu_{t}^{(1/2)}(ds) |\\
&\leq&\int_{0}^{+\infty}T^{\ast}f(x)
|\frac{\partial^{k}}{\partial t^{k}} \mu_{t}^{(1/2)}(ds)|
\leq C_{k} \, T^{\ast}f(x) t^{-k}.
\end{eqnarray*}\ep

\begin{lemma}\label{lematec3}
Given $f\in L^{p}(\gamma_{d}), \alpha\geq 0$ and $k,l$ integers
greater than $\alpha$, then 
$$\|\frac{\partial^{k}P_tf}{\partial
t^{k}}\|_{p,\gamma_d}\leq A_{k}t^{-k+\alpha}\, \mbox{if and only
if} \, \, \|\frac{\partial^{l}P_tf}{\partial
t^{l}}\|_{p,\gamma_d}\leq A_{l}t^{-l+\alpha}.$$ Moreover, if
$A_{k}(f),A_{l}(f)$  are the smallest constants appearing in the
above inequalities then there exist constants $A_{k,l,\alpha}$ and
$D_{k,l,\alpha}$ such that
$$ A_{k,l,\alpha} A_{k}(f)\leq A_{l}(f)\leq
D_{k,l,\alpha}A_{k}(f),$$ for all $f\in L^{p}(\gamma_{d})$.
\end{lemma}

\dem Let us suppose, without lost of generality, that $k\geq l$. We will
prove first the direct implication. For this, we use the
representation of the Poisson-Hermite semigroup (\ref{02}),
$$P_{t}f(x)=\int_{0}^{+\infty}T_{s}f(x)\mu_{t}^{(1/2)}(ds),$$
then differentiating $k$-times with respect to $t$,

$$\frac{\partial^{k}P_{t}f(x)}{\partial
t^{k}}=\int_{0}^{+\infty}T_{s}f(x)
\frac{\partial^{k}}{\partial t^{k}}\mu_{t}^{(1/2)}(ds).$$

Using the identity (\ref{inductdient}) it is easy to prove that for all $m\in
\mathbb{N}$
$$\lim_{t\rightarrow +\infty}\frac{\partial^{m}P_{t}f(x)}{\partial
t^{m}}=0,$$
and therefore given $n\in \mathbb{N}, n>\alpha$
\begin{eqnarray*}
\frac{\partial^{n}P_{t}f(x)}{\partial t^{n}} =-\int_{t}^{+\infty}\frac{\partial^{n+1}P_{s}f(x)}{\partial
s^{n+1}}ds
\end{eqnarray*}
Thus,
\begin{eqnarray*}
\|\frac{\partial^{n}P_tf}{\partial t^{n}}\|_{p,\gamma_d}
&\leq&
\int_{t}^{+\infty}\|\frac{\partial^{n+1}P_sf}{\partial
s^{n+1}}\|_{p,\gamma_d}ds
\leq\int_{t}^{+\infty}A_{n+1}(f)s^{-(n+1)+\alpha}ds\\
&=&\frac{A_{n+1}(f)}{n-\alpha}t^{-n+\alpha}.
\end{eqnarray*}
Then
$$A_{n}(f)\leq\frac{A_{n+1}(f)}{n-\alpha},$$ and as
$n>\alpha$ is arbitrary, then by using the above result $k-l$
times, we get
\begin{eqnarray*}
A_{l}(f)&\leq&
\frac{A_{l+1}(f)}{l-\alpha}\leq\frac{A_{l+2}}{(l-\alpha)(l+1-\alpha)}
\leq...\leq\frac{A_{k}(f)}{(l-\alpha)(l+1-\alpha)...(k-1-\alpha)}\\
&=&D_{k,l,\alpha}A_{k}(f).
\end{eqnarray*}
To prove the converse implication, using again the representation
(\ref{02}), we get,
$$u(x,t_{1}+t_{2})=P_{t_{1}}(P_{t_{2}}f)(x)=\displaystyle\int_{0}^{+\infty}T_{s}(P_{t_{2}}f)(x)\mu_{t_{1}}^{(1/2)}(ds).$$
Therefore, taking  $t=t_{1}+t_{2}$ and differentiating $l$ times
with respect to $t_{2}$ and $k-l$ times with respect to $t_{1}$ we
get
\begin{equation}\label{PoissonkDev}
\frac{\partial^{k}u(x,t)}{\partial
t^{k}}=\int_{0}^{+\infty}T_{s}
(\frac{\partial^{l}P_{t_{2}}f(x)}{\partial
t_{2}^{l}})\frac{\partial^{k-l}}{\partial
t_{1}^{k-l}} \mu_{t_{1}}^{(1/2)}(ds).
\end{equation}

 Thus, using  the inequality (\ref{onesideineq}) and the fact that  the
 Ornstein-Uhlenbeck semigroup is a contraction semigroup, we get
\begin{eqnarray*}
\|\frac{\partial^{k}u(\cdot,t)}{\partial
t^{k}}\|_{p,\gamma_d}&\leq&\int_{0}^{+\infty}\|T_{s}
(\frac{\partial^{l}P_{t_{2}}f}{\partial
t_{2}^{l}})\|_{p,\gamma_d}|\frac{\partial^{k-l}
\mu_{t_{1}}^{(1/2)}}{\partial
t_{1}^{k-l}}(ds)|\leq\|\frac{\partial^{l}P_{t_{2}}f}{\partial
t_{2}^{l}}\|_{p,\gamma_d}\int_{0}^{+\infty}|\frac{\partial^{k-l}}{\partial
t_{1}^{k-l}}\mu_{t_{1}}^{(1/2)}(ds)|\\
&\leq& C_{k-l} \|\frac{\partial^{l}}{\partial
t_{2}^{l}}P_{t_{2}}f\|_{p,\gamma_d}t_{1}^{l-k} \leq C_{k-l}
A_{l}(f)t_{2}^{-l+\alpha}t_{1}^{l-k}.
\end{eqnarray*}
Therefore, taking $t_{1}=t_{2}=\frac{t}{2}$,
$$\displaystyle\|\frac{\partial^{k}u(\cdot,t)}{\partial
t^{k}}\|_{p,\gamma_d}\leq C_{k-l}
A_{l}(f)(\frac{t}{2})^{-k+\alpha},$$ and then,
$$\displaystyle
A_{k}(f)\leq\frac{C_{k-l}}{2^{-k+\alpha}} A_{l}(f).$$\ep

The following technical result will be the key to define Gaussian Besov-Lipschitz spaces,
\begin{lemma}\label{lematec4}
Given $\alpha \geq 0$ and $k,l$ integers greater than $\alpha$. Then
$$\big(\int_{0}^{+\infty}\big(t^{k-\alpha}\|\frac{\partial^{k}P_tf}{\partial
t^{k}}\|_{p,\gamma_d}\big)^{q}\frac{dt}{t}\big)^{\frac{1}{q}}<\infty$$
if and only if
$$\big(\int_{0}^{+\infty}\big(t^{l-\alpha}\|\frac{\partial^{l}P_tf}{\partial
t^{l}}\|_{p,\gamma_d}\big)^{q}\frac{dt}{t}\big)^{\frac{1}{q}}<\infty.
 $$ Moreover,  there exists  constants $A_{k,l,\alpha}, D_{k,l,\alpha}$ such that
\begin{eqnarray*}
D_{k,l,\alpha}\big(\int_{0}^{+\infty}\big(t^{l-\alpha}\|\frac{\partial^{l}P_tf}{\partial
t^{l}}\|_{p,\gamma_d}\big)^{q}\frac{dt}{t}\big)^{\frac{1}{q}}
&\leq&\big(\int_{0}^{+\infty}\big(t^{k-\alpha}\|\frac{\partial^{k}P_tf}{\partial
t^{k}}\|_{p,\gamma_d}\big)^{q}\frac{dt}{t}\big)^{\frac{1}{q}}\\
&\leq&A_{k,l,\alpha}\big(\int_{0}^{+\infty}\big(t^{l-\alpha}\|\frac{\partial^{l}P_tf}{\partial
t^{l}}\|_{p,\gamma_d}\big)^{q}\frac{dt}{t}\big)^{\frac{1}{q}}
\end{eqnarray*}
\end{lemma}
\dem Let us suppose, without lost of generality, that $k\geq l$. We will prove
first the converse implication; from Lemma
 \ref{lematec3}, we have,
\begin{eqnarray*}
\|\frac{\partial^{k}P_tf}{\partial
t^{k}}\|_{p,\gamma_d}&\leq& C_{k-l}\|\frac{\partial^{l}P_{\frac{t}{2}}f}{\partial
(\frac{t}{2})^{l}}\|_{p,\gamma_d}(\frac{t}{2})^{l-k}.
\end{eqnarray*}
Thus,
\begin{eqnarray*}
\big(\int_{0}^{+\infty}\big(t^{k-\alpha}\|\frac{\partial^{k}P_tf}{\partial
t^{k}}\|_{p,\gamma_d}\big)^{q}\frac{dt}{t}\big)^{\frac{1}{q}}&\leq&\frac{C_{k-l}}{2^{l-k}}
\big(\int_{0}^{+\infty}\big(t^{l-\alpha}\|\frac{\partial^{l}P_{t/2}f}{\partial
(\frac{t}{2})^{l}}\|_{p,\gamma_d}\big)^{q}\frac{dt}{t}\big)^{\frac{1}{q}}\\
&=&A_{k,l,\alpha}\big(\int_{0}^{+\infty}\big(s^{l-\alpha}\|\frac{\partial^{l}P_sf}{\partial
s^{l}}\|_{p,\gamma_d}\big)^{q}\frac{ds}{s}\big)^{\frac{1}{q}}
\end{eqnarray*}
with $\displaystyle A_{k,l,\alpha}=C_{k-l} 2^{k-\alpha}$.

For the direct implication, given $n\in \mathbb{N}, \, n>\alpha$,
again using the previous lemma
$$\|\frac{\partial^{n}P_tf}{\partial
t^{n}}\|_{p,\gamma_d}\leq\int_{t}^{+\infty}\|\frac{\partial^{n+1}P_sf}{\partial
s^{n+1}}\|_{p,\gamma_d}ds$$

Therefore, by using the Hardy inequality \cite{se70}
\begin{eqnarray*}
&&\big(\int_{0}^{+\infty}\big(t^{n-\alpha}\|\frac{\partial^{n}P_tf}{\partial
t^{n}}\|_{p,\gamma_d}\big)^{q}\frac{dt}{t}\big)^{\frac{1}{q}}
\quad \quad\quad \quad \quad \quad \quad \quad\quad \quad \quad
\quad\\
&&\quad \quad\quad \quad \quad \quad \quad \leq\big(\int_{0}^{+\infty}\big(t^{n-\alpha}\int_{t}^{+\infty}\|\frac{\partial^{n+1}P_sf}{\partial
s^{n+1}}\|_{p,\gamma_d}ds\big)^{q}\frac{dt}{t}\big)^{\frac{1}{q}}\\
&&\quad \quad\quad \quad \quad \quad \quad=\big(\int_{0}^{+\infty}\big(\int_{t}^{+\infty}\|\frac{\partial^{n+1}P_sf}{\partial
s^{n+1}}\|_{p,\gamma_d}ds\big)^{q}t^{(n-\alpha)q-1}dt\big)^{\frac{1}{q}}\\
&& \quad \quad\quad \quad \quad \quad \quad \leq\frac{1}{n-\alpha}\big(\int_{0}^{+\infty}\big(s^{n+1-\alpha}\|\frac{\partial^{n+1}P_sf}{\partial
s^{n+1}}\|_{p,\gamma_d}\big)^{q}\frac{ds}{s}\big)^{\frac{1}{q}}.
\end{eqnarray*}
Now, as $n>\alpha$ is arbitrary,  using the above result $k-l$
times

\begin{eqnarray*}
&&\big(\int_{0}^{+\infty}\big(t^{l-\alpha}\|\frac{\partial^{l}P_tf}{\partial
t^{l}}\|_{p,\gamma_d}\big)^{q}\frac{dt}{t}\big)^{\frac{1}{q}}
\quad \quad\quad \quad \quad \quad \quad \quad\quad \quad \quad
\quad \\
&& \quad \quad\quad \quad \quad \quad \quad\leq\frac{1}{l-\alpha}\big(\int_{0}^{+\infty}\big(t^{l+1-\alpha}\|\frac{\partial^{l+1}P_tf}{\partial
t^{l+1}}\|_{p,\gamma_d}\big)^{q}\frac{dt}{t}\big)^{\frac{1}{q}}\\
&&\quad \quad\quad \quad \quad \quad \quad \leq \frac{1}{(l-\alpha).(l+1-\alpha)}\big(\int_{0}^{+\infty}\big(t^{l+2-\alpha}\|\frac{\partial^{l+2}P_tf}{\partial t^{l+2}}\|_{p,\gamma_d}\big)^{q}\frac{dt}{t}\big)^{\frac{1}{q}}\\
&&\quad \quad\quad \quad \quad \quad \quad ...\\
&& \quad \quad\quad \quad \quad \quad \quad \leq D_{k,l,\alpha}\big(\int_{0}^{+\infty}\big(t^{k-\alpha}\|
\frac{\partial^{k}P_tf}{\partial
t^{k}}\|_{p,\gamma_d}\big)^{q}\frac{dt}{t}\big)^{\frac{1}{q}}
\end{eqnarray*}
where $\displaystyle
D_{k,l,\alpha}=\frac{1}{(l-\alpha).(l+1-\alpha)...(k-1-\alpha)}$. \ep \\

Now, following the classical case, see for instance \cite{fm91},
\cite{se70}, \cite{trie1} and \cite{trie2}, we  are going to
define the Gaussian Besov-Lipschitz $B_{p,q}^{\alpha}(\gamma_d)$
spaces or Besov-Lipschitz spaces  for Hermite polynomial
expansions,

\begin{defi}
Let $\alpha \geq 0$, $k$ be the smallest integer greater than
$\alpha$, and $1\leq p,q \leq\infty$. For $1 \leq q < \infty$ the Gaussian Besov-Lipschitz
space $B_{p,q}^{\alpha}(\gamma_d)$ are defined as the set of functions  $f \in
L^p(\gamma_d)$ for which
\begin{equation}\label{e15}
 \left( \int_0^{\infty} (t^{k-\alpha} \left\|
\frac{\partial^{k}P_t f}{\partial t^{k}} \right\|_{p,\gamma_d}
)^{q} \frac{dt}{t} \right) ^{1/q}  < \infty.
\end{equation}
The norm of $f \in B_{p,q}^{\alpha}(\gamma_d)$ is defined as
\begin{equation}
\left\| f \right\|_{B_{p,q}^{\alpha}}: =  \left\| f \right\|_{p,
\gamma_d} +\left( \int_0^{\infty} (t^{k-\alpha} \left\|
\frac{\partial^{k}   P_t f}{\partial t^{k}} \right\|_{p,\gamma_d}
)^{q} \frac{dt}{t} \right) ^{1/q}
\end{equation}
 For $q=\infty$ the Gaussian Besov-Lipschitz space
$B_{p,\infty}^{\alpha}(\gamma_d)$ are defined  as the set of  functions  $f \in
L^p(\gamma_d)$ for which exists a constant $A$ such that\\
$$\|\frac{\partial^{k}P_t f}{\partial
t^{k}}\|_{p,\gamma_d}\leq At^{-k+\alpha}$$ and then the norm of $f
\in B_{p,\infty}^{\alpha}(\gamma_d)$ is defined as
\begin{equation}
\left\| f \right\|_{B_{p,\infty}^{\alpha}}: =  \left\| f
\right\|_{p, \gamma_d} +A_{k}(f),
\end{equation}
where $A_{k}(f)$ is the smallest constant $A$ appearing in the
above inequality.

In particular, the space $B_{\infty,\infty}^{\alpha}(\gamma_d)$ is the Gaussian Lipschitz space $Lip_{\alpha}(\gamma_d)$.
\end{defi}

Lemma \ref{lematec4} show us that we could
have replaced $k$ by any other integer $l$ greater than $\alpha$ and the
resulting norms are equivalent.  

In what follows, we need the following
technical result about $L^p(\gamma_d)$-norms of the derivatives of
the Poisson-Hermite semigroup,

\begin{lemma}\label{kdecay} 
Suppose $f\in L^{p}(\gamma_{d})$, then for any integer $k$ the function
$\displaystyle\|\frac{\partial^{k}P_t f}{\partial
t^{k}}\|_{p,\gamma_{d}}$ is a non-increasing function of $t$, for
$0<t<+\infty$. Moreover,
\begin{equation}\label{kdecayine}
\|\frac{\partial^{k}P_t f}{\partial t^{k}}\|_{p,\gamma_d}\leq
C \| f\|_{p,\gamma_d}t^{-k}, t>0
\end{equation}

\end{lemma}
\dem Let us consider first the case $k=0$. Let us fix $t_{1},t_{2}>0$,
by using the semigroup property we get

$$
u(x,t_{1}+t_{2})=P_{t_{1}+t_{2}}f(x) =P_{t_{1}}(P_{t_{2}}f(x)) =
P_{t_{1}}(u(x,t_{2}))
$$
Therefore, by definition of $P_{t}$, Jensen's inequality and the
invariance of $\gamma_{d}$
\begin{eqnarray*}
\int_{{\mathbb{R}}^{d}}|u(x,t_{1}+t_{2})|^{p}\gamma_{d}(dx)&=&\int_{{\mathbb{R}}^{d}}|\int_{{\mathbb{R}}^{d}}p(t_{1},x,y) u(y,t_{2})dy|^{p}\gamma_{d}(dx)\\
&\leq&\int_{{\mathbb{R}}^{d}}\big(\int_{{\mathbb{R}}^{d}}p(t_{1},x,y)|u(y,t_{2})|^{p}dy\big)\gamma_{d}(dx)\\
&=&\int_{{\mathbb{R}}^{d}}P_{t_{1}}(|u(x,t_{2})|^{p})\gamma_{d}(dx)=\int_{{\mathbb{R}}^{d}}|u(x,t_{2})|^{p}\gamma_{d}(dx).
\end{eqnarray*}
Thus
$$\|P_{t_{1}+t_{2}}f\|_{p,\gamma_{d}}\leq\|P_{t_{2}}f\|_{p,\gamma_{d}}.$$

Now to prove the general case, $k >0$. Differentiating the identity
 $u(x,t_{1}+t_{2})=P_{t_{1}}(u(x,t_{2}))$ $k$-times with respect to $t_{2}$ to get

$$ \frac{\partial^{k}u(x,t_{1}+t_{2})}{\partial (t_{1}+t_{2})^{k}}=P_{t_{1}}(\frac{\partial^{k}u(x,t_{2})}{\partial
t_{2}^{k}})$$ and then use a analogous argument to the one above.

In order to prove (\ref{kdecayine}) we use again the
representation (\ref{02}) of the Poisson-Hermite  semigroup and
differentiating  it $k$-times with respect to $t$ we get

$$\frac{\partial^{k}P_tf(x)}{\partial
t^{k}}=\int_{0}^{+\infty}T_{s}f(x)\frac{\partial^{k}}{\partial
t^{k}}\mu_{t}^{(1/2)}(ds),$$
 thus, by Minkowski's integral inequality, the contractive property of the
  Ornstein-Uhlenbeck semigroup and inequality (\ref{onesideineq}), we get for $t>0$
\begin{eqnarray*}
\|\frac{\partial^{k}P_tf}{\partial
t^{k}}\|_{p,\gamma_d}&\leq&\int_{0}^{+\infty}\|T_{s}f
\frac{\partial^{k}}{\partial
t^{k}}\mu_{t}^{(1/2)}(ds)\|_{p,\gamma_d}=\int_{0}^{+\infty}\|T_{s}f\|_{p,\gamma_d}|\frac{\partial^{k}}{\partial
t^{k}}\mu_{t}^{(1/2)}(ds)|\\
&\leq&\|f\|_{p,\gamma_d}\int_{0}^{+\infty}|\frac{\partial^{k}
}{\partial t^{k}}\mu_{t}^{(1/2)}(ds)| \leq
\frac{C_{k}}{t^{k}} \|f\|_{p,\gamma_d}.
\end{eqnarray*} \ep \\

Let us study some inclusions among the Gaussian Besov-Lipschitz spaces,

\begin{proposition}\label{incluBesov}
The inclusion $B_{p,q_1}^{\alpha_{1}}(\gamma_d)\subset
B_{p,q_2}^{\alpha_{2}}(\gamma_d)$ holds if either:
\begin{enumerate}
\item [i)]   $\alpha_{1}>\alpha_{2}>0$ ($q_{1}$ and $q_{2}$ need
not to be related), or \item[ii)] If $\alpha_{1}=\alpha_{2}$ and
$q_{1}\leq q_{2}$
\end{enumerate}
\end{proposition}
\dem In order to prove  $ii)$, we set
$A=\big(\displaystyle\int_{0}^{+\infty}\big(t^{k-\alpha}\|\frac{\partial^{k}P_t
f}{\partial
t^{k} }\|_{p,\gamma_d}\big)^{q_{1}}\frac{dt}{t}\big)^{\frac{1}{q_{1}}}$\\
Now, fixing $t_{0}>0$
$$\displaystyle\int_{\frac{t_{0}}{2}}^{t_{0}}\big(t^{k-\alpha}\|\frac{\partial^{k}P_t f}{\partial
t^{k} }\|_{p,\gamma_d}\big)^{q_{1}}\frac{dt}{t}\leq A^{q_{1}}.$$
By Lemma \ref{kdecay},
$\displaystyle\|\frac{\partial^{k}P_t f}{\partial t^{k}
}\|_{p,\gamma_d}$ takes its minimum value at the upper end point
$(t=t_{0})$ of the above integral . So we get
$$\displaystyle\|\frac{\partial^{k}P_{t_{0}} f}{\partial t^{k}
}\|_{p,\gamma_d}^{q_{1}}\displaystyle\int_{\frac{t_{0}}{2}}^{t_{0}}t^{(k-\alpha)q_{1}}\frac{dt}{t}\leq
A^{q_{1}}.$$ That is $\displaystyle\|\frac{\partial^{k}P_{t_{0}}
f}{\partial t^{k}}\|_{p,\gamma_d}\leq CAt_{0}^{-k+\alpha}$ but
since $t_{0}$ is arbitrary then
$$\displaystyle\|\frac{\partial^{k}P_t f}{\partial t^{k}
}\|_{p,\gamma_d}\leq CAt^{-k+\alpha},$$
 for all $t>0$.
In other words $f\in B_{p,q_1}^{\alpha}$ implies also that $f\in
B_{p,\infty}^{\alpha}$. Thus, as $q_{2}\geq q_{1}$

\begin{eqnarray*}
\int_{0}^{+\infty}\big(t^{k-\alpha}\|\frac{\partial^{k}P_t
f}{\partial t^{k}
}\|_{p,\gamma_d}\big)^{q_{2}}\frac{dt}{t}&=&\int_{0}^{+\infty}\big(t^{k-\alpha}\|\frac{\partial^{k}P_t
f}{\partial t^{k}
}\|_{p,\gamma_d}\big)^{q_{2}-q_{1}}\big(t^{k-\alpha}\|\frac{\partial^{k}P_t
f}{\partial
t^{k} }\|_{p,\gamma_d}\big)^{q_{1}}\frac{dt}{t}\\
&\leq&(CA)^{q_{2}-q_{1}}\int_{0}^{+\infty}\big(t^{k-\alpha}\|\frac{\partial^{k}P_t
f}{\partial
t^{k} }\|_{p,\gamma_d}\big)^{q_{1}}\frac{dt}{t}\\
&=&(CA)^{q_{2}-q_{1}} A^{q_{1}}=CA^{q_{2}}<+\infty,
\end{eqnarray*}
and therefore $f\in B_{p,q_2}^{\alpha}$.\\

Now in order to prove part $i)$, \ by Lemma \ref{kdecay} we have
$$\|\frac{\partial^{k}P_t f}{\partial t^{k}}\|_{p,\gamma_d}\leq
Ct^{-k}, \, t>0.$$ Now given $f\in B_{p,q_1}^{\alpha_{1}}$, taking
again
$$A=\big(\displaystyle\int_{0}^{+\infty}\big(t^{k-\alpha_{1}}\|\frac{\partial^{k}P_t f}{\partial
t^{k}
}\|_{p,\gamma_d}\big)^{q_{1}}\frac{dt}{t}\big)^{\frac{1}{q_{1}}},$$
we get as in part $ii)$
$$\displaystyle\|\frac{\partial^{k}P_t f}{\partial t^{k}
}\|_{p,\gamma_d}\leq CAt^{-k+\alpha_{1}},$$ for all $t>0$. Now,
\begin{eqnarray*}
\int_{0}^{+\infty}\big(t^{k-\alpha_{2}}\|\frac{\partial^{k}P_t
f}{\partial t^{k}
}\|_{p,\gamma_d}\big)^{q_{2}}\frac{dt}{t}&=&\int_{0}^{1}\big(t^{k-\alpha_{2}}\|\frac{\partial^{k}P_t
f}{\partial t^{k}
}\|_{p,\gamma_d}\big)^{q_{2}}\frac{dt}{t}+\int_{1}^{+\infty}\big(t^{k-\alpha_{2}}\|\frac{\partial^{k}P_t
f}{\partial
t^{k} }\|_{p,\gamma_d}\big)^{q_{2}}\frac{dt}{t}\\
&=&I+II.
\end{eqnarray*}
Now,
\begin{eqnarray*}
I&=&\int_{0}^{1}t^{(k-\alpha_{2})q_{2}}\|\frac{\partial^{k}P_t
f}{\partial t^{k} }\|_{p,\gamma_d}^{q_{2}}\frac{dt}{t}
\leq \int_{0}^{1}t^{(k-\alpha_{2})q_{2}}(CA)^{q_{2}}t^{(\alpha_{1}-k)q_{2}}\frac{dt}{t}\\
&=&(CA)^{q_{2}}\int_{0}^{1}t^{(\alpha_{1}-\alpha_{2})q_{2}}\frac{dt}{t}=CA^{q_{2}},
\end{eqnarray*}
and
\begin{eqnarray*}
II&=&\int_{1}^{+\infty}t^{(k-\alpha_{2})q_{2}}\|\frac{\partial^{k}P_t
f}{\partial t^{k} }\|_{p,\gamma_d}^{q_{2}}\frac{dt}{t}
\leq\int_{1}^{+\infty}t^{(k-\alpha_{2})q_{2}}C^{q_{2}}t^{-kq_{2}}\frac{dt}{t}\\
&=&C^{q_{2}}\int_{1}^{+\infty}t^{-\alpha_{2}q_{2}}\frac{dt}{t}=C.
\end{eqnarray*}
Hence,
\begin{eqnarray*}
\int_{0}^{+\infty}\big(t^{k-\alpha_{2}}\|\frac{\partial^{k}P_t
f}{\partial t^{k}
}\|_{p,\gamma_d}\big)^{q_{2}}\frac{dt}{t}&<&+\infty,
\end{eqnarray*}
and so $f\in B_{p,q_2}^{\alpha_{2}}$. \ep\\

The following technical result will be the key to define Gaussian Triebel-Lizorkin spaces,
\begin{lemma}\label{TLind}
Let $\alpha \geq 0$ and $k,l$ integers such that $k\geq l>\alpha$. Then
$$
\|\big(\int_{0}^{+\infty}\big(t^{k-\alpha}|\frac{\partial^{k}}{\partial
t^{k}}P_{t}f|\big)^{q}\frac{dt}{t}\big)^{\frac{1}{q}}\|_{p,\gamma}<\infty $$
if and only if
$$\|\big(\int_{0}^{+\infty}\big(t^{l-\alpha}|\frac{\partial^{l}}{\partial
t^{l}}P_{t}f|\big)^{q}\frac{dt}{t}\big)^{\frac{1}{q}}\|_{p,\gamma}<\infty.
$$
Moreover, there exists  constants $A_{k,l,\alpha}, D_{k,l,\alpha}$ such that
\begin{eqnarray*}
D_{k,l,\alpha} \|\big(\int_{0}^{+\infty}\big(t^{l-\alpha}|\frac{\partial^{l}}{\partial
t^{l}}P_{t}f|\big)^{q}\frac{dt}{t}\big)^{\frac{1}{q}}\|_{p,\gamma}&\leq&\|\big(\int_{0}^{+\infty}\big(t^{k-\alpha}|
\frac{\partial^{k}}{\partial
t^{k}}P_{t}f|\big)^{q}\frac{dt}{t}\big)^{\frac{1}{q}}\|_{p,\gamma} \\
&\leq& A_{k,l,\alpha}  \|\big(\int_{0}^{+\infty}\big(t^{l-\alpha}|\frac{\partial^{l}}{\partial
t^{l}}P_{t}f|\big)^{q}\frac{dt}{t}\big)^{\frac{1}{q}}\|_{p,\gamma}.
\end{eqnarray*}  
\end{lemma}
\dem

 Let $n\in \mathbb{N}$ such that $n>\alpha$. Then it can be proved that
$$|\frac{\partial^{n}}{\partial
t^{n}}P_{t}f(x)|\leq\int_{t}^{+\infty}|\frac{\partial^{n+1}}{\partial
s^{n+1}}P_{s}f(x)|ds$$

Then by Hardy's inequality,
\begin{eqnarray*}
\big(\int_{0}^{+\infty}\big(t^{n-\alpha}|\frac{\partial^{n}}{\partial
t^{n}}P_{t}f(x)|\big)^{q}\frac{dt}{t}\big)^{\frac{1}{q}}
&\leq&\big(\int_{0}^{+\infty}\big(t^{n-\alpha}\int_{t}^{+\infty}|\frac{\partial^{n+1}}{\partial
s^{n+1}}P_{s}f(x)|ds\big)^{q}\frac{dt}{t}\big)^{\frac{1}{q}}\\
&\leq&\frac{1}{n-\alpha}\big(\int_{0}^{+\infty}\big(s|\frac{\partial^{n+1}}{\partial
s^{n+1}}P_{s}f(x)|\big)^{q}s^{(n-\alpha)q-1}ds\big)^{\frac{1}{q}}\\
&=&\frac{1}{n-\alpha}\big(\int_{0}^{+\infty}\big(s^{n+1-\alpha}|\frac{\partial^{n+1}}{\partial
s^{n+1}}P_{s}f(x)|\big)^{q}\frac{ds}{s}\big)^{\frac{1}{q}}.
\end{eqnarray*}
Now as  $n>\alpha$ is arbitrary, iterating the previous argument  $k-l$ times, we have
\begin{eqnarray*}
\big(\int_{0}^{+\infty}\big(t^{l-\alpha}|\frac{\partial^{l}}{\partial
t^{l}}P_{t}f(x)|\big)^{q}\frac{dt}{t}\big)^{\frac{1}{q}}
&\leq&\frac{1}{l-\alpha}\big(\int_{0}^{+\infty}\big(t^{l+1-\alpha}|\frac{\partial^{l+1}}{\partial
t^{l+1}}P_{t}f(x)|\big)^{q}\frac{dt}{t}\big)^{\frac{1}{q}}\\
&\leq&\frac{1}{(l-\alpha).(l+1-\alpha)}\big(\int_{0}^{+\infty}\big(t^{l+2-\alpha}|\frac{\partial^{l+2}}{\partial
t^{l+2}}P_{t}f(x)|\big)^{q}\frac{dt}{t}\big)^{\frac{1}{q}}\\
&&...\\
&\leq&C_{k,l,\alpha}\big(\int_{0}^{+\infty}\big(t^{k-\alpha}|
\frac{\partial^{k}}{\partial
t^{k}}P_{t}f(x)|\big)^{q}\frac{dt}{t}\big)^{\frac{1}{q}}
\end{eqnarray*}
where $\displaystyle C_{k,l,\alpha}=\frac{1}{(l-\alpha)(l+1-\alpha)...(k-1-\alpha)}.$
Thus
\begin{eqnarray*}
D_{k,l,\alpha}\|\big(\int_{0}^{+\infty}\big(t^{l-\alpha}|\frac{\partial^{l}}{\partial
t^{l}}P_{t}f|\big)^{q}\frac{dt}{t}\big)^{\frac{1}{q}}\|_{p,\gamma}&\leq&\|\big(\int_{0}^{+\infty}\big(t^{k-\alpha}|
\frac{\partial^{k}}{\partial
t^{k}}P_{t}f|\big)^{q}\frac{dt}{t}\big)^{\frac{1}{q}}\|_{p,\gamma},
\end{eqnarray*}
where $D_{k,l,\alpha} = 1/C_{k,l,\alpha}.$

The converse inequality is also obtained by an inductive argument from the case $k=l+1$. Let us remember (\ref{PoissonkDev}),
$$\frac{\partial^{k}u(x,t)}{\partial
t^{k}}=\int_{0}^{+\infty}T_{s}
(\frac{\partial^{l}P_{t_{2}}f(x)}{\partial
t_{2}^{l}})\frac{\partial^{k-l}}{\partial
t_{1}^{k-l}} \mu_{t_{1}}^{(1/2)}(ds),$$
and since, from (\ref{inductdient}), $\displaystyle\frac{\partial }{\partial
t_{1}}\mu_{t_{1}}^{(1/2)} (ds)=\big(t_{1}^{-1}-\frac{t_{1}}{2s}\big)\mu_{t_{1}}^{(1/2)}(ds)$
we get
\begin{eqnarray*}
|\displaystyle\frac{\partial^{k}u(x,t)}{\partial
t^{k}}|&\leq&\int_{0}^{+\infty}T_{s}
(|\frac{\partial^{l}P_{t_{2}}f(x)}{\partial
t_{2}^{l}}|) |\big(t_{1}^{-1}-\frac{t_{1}}{2s}\big)|\mu_{t_{1}}^{(1/2)}(ds)\\
&\leq&t_{1}^{-1}\int_{0}^{+\infty}T_{s}
(|\frac{\partial^{l}P_{t_{2}}f(x)}{\partial
t_{2}^{l}}|)\mu_{t_{1}}^{(1/2)}(ds)+\frac{t_{1}}{2}\int_{0}^{+\infty}T_{s}
(|\frac{\partial^{l}P_{t_{2}}f(x)}{\partial
t_{2}^{l}}|)\frac{1}{s} \mu_{t_{1}}^{(1/2)}(ds).
\end{eqnarray*}
Therefore
\begin{eqnarray*}
\big(\int_{0}^{+\infty}\big(t_{2}^{k-\alpha}|\displaystyle\frac{\partial^{k}u(x,t)}{\partial
t^{k}}|\big)^{q}\frac{dt_{2}}{t_{2}}\big)^{1/q}&\leq&C_{q}\big[\big(\int_{0}^{+\infty}\big(t_{2}^{k-\alpha}t_{1}^{-1} \int_{0}^{+\infty}T_{s}
(|\frac{\partial^{l}P_{t_{2}}f(x)}{\partial
t_{2}^{l}}|)\mu_{t_{1}}^{(1/2)}(ds)\big)^{q}\frac{dt_{2}}{t_{2}}\big)^{1/q}\\
&+&\big(\int_{0}^{+\infty}\big(t_{2}^{k-\alpha}\frac{t_{1}}{2}\int_{0}^{+\infty}T_{s}
(|\frac{\partial^{l}P_{t_{2}}f(x)}{\partial
t_{2}^{l}}|)\frac{1}{s} \mu_{t_{1}}^{(1/2)}(ds)\big)^{q}\frac{dt_{2}}{t_{2}}\big)^{1/q}\big]\\
&=&I + II
\end{eqnarray*}
Now using twice Minkowski integral inequality (since $T_s$ is an integral transformation with positive kernel) and the fact that $\mu_{t_{1}}^{(1/2)}(ds)$ is a probability, we get
\begin{eqnarray*}
I &=&C_{q} \big(\int_{0}^{+\infty}\big(t_{2}^{k-\alpha}t_{1}^{-1}\big)^{q}\big(\int_{0}^{+\infty}T_{s}
(|\frac{\partial^{l}P_{t_{2}}f(x)}{\partial t_{2}^{l}}|)\mu_{t_{1}}^{(1/2)}(ds)\big)^{q}\frac{dt_{2}}{t_{2}}\big)^{1/q}\\
&\leq& C_{q}\int_{0}^{+\infty}\big(\int_{0}^{+\infty}\big(t_{2}^{k-\alpha}t_{1}^{-1}\big)^{q}\big(T_{s}
(|\frac{\partial^{l}P_{t_{2}}f(x)}{\partial
t_{2}^{l}}|)\big)^{q}\frac{dt_{2}}{t_{2}}\big)^{1/q}\mu_{t_{1}}^{(1/2)}(ds)\\
&\leq&C_{q}\int_{0}^{+\infty}T_{s}\big(\big(\int_{0}^{+\infty}\big(t_{2}^{k-\alpha}t_{1}^{-1}\big)^{q}\big(|\frac{\partial^{l}P_{t_{2}}f(x)}{\partial
t_{2}^{l}}|\big)^{q}\frac{dt_{2}}{t_{2}}\big)^{1/q}\big)\mu_{t_{1}}^{(1/2)}(ds)\\
&\leq&C_{q} T^{\ast}\big(\big(\int_{0}^{+\infty}\big(t_{2}^{k-\alpha}t_{1}^{-1}\big)^{q}.\big(|\frac{\partial^{l}P_{t_{2}}f(x)}{\partial
t_{2}^{l}}|\big)^{q}\frac{dt_{2}}{t_{2}}\big)^{1/q}\big)
\end{eqnarray*}
and using the same argument for $(II)$ and (\ref{intoneside}), we have
\begin{eqnarray*}
II &\leq &C_qT^{\ast}\big(\big(\int_{0}^{+\infty}\big(t_{2}^{k-\alpha}t_{1}\big)^{q}\big(|\frac{\partial^{l}P_{t_{2}}f(x)}{\partial
t_{2}^{l}}|\big)^{q}\frac{dt_{2}}{t_{2}}\big)^{1/q}\big)\frac{1}{t_{1}^{2}}\\
&=&C_q T^{\ast}\big(\big(\int_{0}^{+\infty}\big(t_{2}^{k-\alpha}t_{1}^{-1}\big)^{q}\big(|\frac{\partial^{l}P_{t_{2}}f(x)}{\partial t_{2}^{l}}|\big)^{q}\frac{dt_{2}}{t_{2}}\big)^{1/q}\big).
\end{eqnarray*}
Taking  $t_{1}=t_{2}=\frac{t}{2}$ and changing the variable, we get
\begin{eqnarray*}
 I &\leq&C_q T^{\ast}\big(\big(\int_{0}^{+\infty}\big(t^{l-\alpha}\big)^{q}\big(|\frac{\partial^{l}P_{t}f(x)}{\partial t^{l}}|\big)^{q}\frac{dt}{t}\big)^{1/q}\big)
\end{eqnarray*}
and
\begin{eqnarray*}
II&\leq&C_q T^{\ast}\big(\big(\int_{0}^{+\infty}\big(t^{l-\alpha}\big)^{q}\big(|\frac{\partial^{l}P_{t}f(x)}{\partial t^{l}}|\big)^{q}\frac{dt}{t}\big)^{1/q}\big).
\end{eqnarray*}
Hence, by the $L^p$ boundedness of $T^*$
\begin{eqnarray*}
\|\big(\int_{0}^{+\infty}\big(t^{k-\alpha}|\displaystyle\frac{\partial^{k}u(x,t)}{\partial
t^{k}}|\big)^{q}\frac{dt}{t}\big)^{1/q}\|_{p,\gamma}&\leq&C_{q,k,\alpha} \|T^{\ast}\big(\big(\int_{0}^{+\infty}\big(t^{l-\alpha}|\frac{\partial^{l}P_{t}f(x)}{\partial
u^{l}}|\big)^{q}\frac{dt}{t}\big)^{1/q}\big)\|_{p,\gamma}\\
&+&C_q\|T^{\ast}\big(\big(\int_{0}^{+\infty}\big(t^{l-\alpha}|\frac{\partial^{l}P_{t}f(x)}{\partial
u^{l}}|\big)^{q}\frac{dt}{t}\big)^{1/q}\big)\|_{p,\gamma}\big)\\
&\leq&C_{k,\alpha,q}\|\big(\int_{0}^{+\infty}\big(t^{l-\alpha}|\frac{\partial^{l}P_{t}f(x)}{\partial
t^{l}}|\big)^{q}\frac{dt}{t}\big)^{1/q}\|_{p,\gamma}.
\end{eqnarray*}
 \ep

Now, we can  introduce the Gaussian Triebel-Lizorkin
 spaces  $F_{p,q}^{\alpha}(\gamma_{d})$ following the classical case (see
\cite{fm91}, \cite{trie1} and \cite{trie2}),
\begin{defi}
Let $\alpha \geq  0$, $k$ be the smallest integer
greater than $\alpha$, and $1\leq p,q<\infty$.  The Gaussian  Triebel-Lizorkin
 space $F_{p,q}^{\alpha}(\gamma_{d})$  is the set of functions $f\in
L^{p}(\gamma_{d})$ for which
\begin{equation}\label{e16}
\left\| \left( \int_0^{\infty} (t^{k-\alpha}
\left|\frac{\partial^{k}P_t f}{\partial t^{k}} \right|) ^{q}
\frac{dt}{t} \right) ^{1/q} \right\|_{p,\gamma_d}<\infty.
\end{equation}
The norm of $f \in F_{p,q}^{\alpha}(\gamma_d)$ is defined as
\begin{equation}
\left\| f \right\|_{F_{p,q}^{\alpha}}: =  \left\| f \right\|_{p,
\gamma_d} +\left\| \left( \int_0^{\infty} (t^{k-\alpha}
\left|\frac{\partial^{k}P_t f}{\partial t^{k}} \right|) ^{q}
\frac{dt}{t} \right) ^{1/q} \right\|_{p,\gamma_d}.
\end{equation}
\end{defi}
Observe that by Lemma \ref{TLind} the definition of $F_{p}^{\alpha,q}(\gamma_d)$ does not depend on which  $k>\alpha$ is chosen and the resulting norms are equivalent.

In \cite{iris} the notion of homogeneous Gaussian Besov-Lipschitz
and homogeneous Gaussian Triebel-Lizorkin spaces were considered.
Nevertheless the definitions of those spaces given there appear to
be wrong in the case that  $\alpha > 1$. On the other hand,  J.
Epperson  \cite{ep} has considered Triebel-Lizorkin spaces with
respect to the Hermite functions expansions which are different to the spaces
that we are considering in this article related to  Hermite
polynomial expansions.

Let us observe that by the $L^p(\gamma_d)$-continuity of the Gaussian
Littlewood-Paley  $g_1$- function, see \cite{gu}
\begin{equation}
g_1(f)(x) = \left( \int_0^{\infty} t \left|\frac{\partial P_t f}{\partial t} \right|^2 dt\right)^{1/2}
\end{equation}
it is inmediate to see that for $1< p<\infty$
$$ L^p(\gamma_d) = F^0_{p,2}(\gamma_d),$$
and by the trivial identification of the $L^p$ spaces with the Hardy spaces, see \cite{forur}, we have also
$$ H^p(\gamma_d) = F^0_{p,2}(\gamma_d),$$

For Gaussian Triebel-Lizorkin spaces we have the following inclusion result, which is analogous to Proposition \ref{incluBesov} i),

\begin{proposition}\label{incluTriebel}
The inclusion $F_{p,q_1}^{\alpha_{1}}(\gamma_d)\subset
F_{p,q_2}^{\alpha_{2}}(\gamma_d)$ holds for $\alpha_1 > \alpha_2>0$ and
$q_1 \geq q_2$.
\end{proposition}
\dem
 Let us consider $f\in F_{p}^{\alpha_{1},q_{1} }(\gamma_d)$. Then
\begin{eqnarray*}
&& \big(\int_{0}^{+\infty}\big(t^{k-\alpha_{2}}|\frac{\partial^{k}P_{t}f(x)}{\partial
t^{k}}|\big)^{q_{2}}\frac{dt}{t}\big)^{\frac{1}{q_{2}}}\quad \quad\quad \quad \quad \quad
\quad\quad \quad\quad \quad\\
&&\quad \quad \quad \quad \quad = \big(\int_{0}^{1} \big(t^{k-\alpha_{2}}|\frac{\partial^{k}P_{t}f(x)}{\partial
t^{k}}|\big)^{q_{2}}\frac{dt}{t}+\int_{1}^{+\infty}\big(t^{k-\alpha_{2}}|\frac{\partial^{k}P_{t}f(x)}{\partial
t^{k}}|\big)^{q_{2}}\frac{dt}{t}\big)^{\frac{1}{q_{2}}}\\
&&\quad \quad \quad \quad \quad \leq\big(\int_{0}^{1} \big(t^{k-\alpha_{2}}|\frac{\partial^{k}P_{t}f(x)}{\partial
t^{k}}|\big)^{q_{2}}\frac{dt}{t}\big)^{\frac{1}{q_{2}}}+ \big(\int_{1}^{+\infty}\big(t^{k-\alpha_{2}}|\frac{\partial^{k}P_{t}f(x)}{\partial
t^{k}}|\big)^{q_{2}}\frac{dt}{t}\big)^{\frac{1}{q_{2}}}\\
&&\quad \quad \quad \quad \quad =  I + II.
\end{eqnarray*}

Let us observe that for the first term $I$, the case $q_{1}=q_{2}$ is immediate since as $t<1$ ,  $t^{k-\alpha_{2}} < t^{k-\alpha_{1}}$ and then
$$I^{q_{2}} \leq  \int_{0}^{+\infty}\big(t^{k-\alpha_{1}}|\frac{\partial^{k}P_{t}f(x)}{\partial
t^{k}}|\big)^{q_{1}}\frac{dt}{t}.$$
 Now, in the case $q_{1}>q_{2}$ taking  $r=\displaystyle\frac{q_{1}}{q_{2}}$,
$s=\displaystyle\frac{q_{1}}{q_{1}-q_{2}}$ then $r,s> 1$ and
$\displaystyle\frac{1}{r}+\frac{1}{s}=1$, then, by Holder's inequality
\begin{eqnarray*}
I^{q_{2}}&=&\int_{0}^{1}t^{(\alpha_{1}-\alpha_{2})q_{2}}\big(t^{k-\alpha_{1}}|\frac{\partial^{k}P_{t}f(x)}{\partial
t^{k}}|\big)^{q_{2}}\frac{dt}{t}\leq \big(\int_{0}^{1}t^{(\alpha_{1}-\alpha_{2})q_{2}s}\frac{dt}{t}\big)^{\frac{1}{s}} \big(\int_{0}^{1}\big(t^{k-\alpha_{1}}|\frac{\partial^{k}P_{t}f(x)}{\partial
t^{k}}|\big)^{q_{2}r}\frac{dt}{t}\big)^{\frac{1}{r}}\\
&=&\frac{1}{(\alpha_{1}-\alpha_{2})q_{2}s}\big(\int_{0}^{1}\big(t^{k-\alpha_{1}}|\frac{\partial^{k}P_{t}f(x)}{\partial
t^{k}}|\big)^{q_{1}}\frac{dt}{t}\big)^{\frac{q_{2}}{q_{1}}} \leq C \big(\int_{0}^{+\infty}\big(t^{k-\alpha_{1}}|\frac{\partial^{k}P_{t}f(x)}{\partial
t^{k}}|\big)^{q_{1}}\frac{dt}{t}\big)^{\frac{q_{2}}{q_{1}}}.
\end{eqnarray*}
Now for the second term $II$, using Lemma \ref{lematec5}, we have
\begin{eqnarray*}
II&=&\big(\int_{1}^{+\infty}\big(t^{k-\alpha_{2}}|\frac{\partial^{k}P_{t}f(x)}{\partial
t^{k}}|\big)^{q_{2}}\frac{dt}{t}\big)^{\frac{1}{q_{2}}}\leq C\,T^{\ast}f(x) \big(\int_{1}^{+\infty}\big(t^{k-\alpha_{2}}t^{-k}\big)^{q_{2}}\frac{dt}{t}\big)^{\frac{1}{q_{2}}}\\
&=&C \, T^{\ast}f(x) \big(\int_{1}^{+\infty}t^{-\alpha_{2} q_{2}}\frac{dt}{t}\big)^{\frac{1}{q_{2}}}=  C\,T^{\ast}f(x).\end{eqnarray*}
Then, using the $L^{p}(\gamma_{d})$
continuity of $T^{\ast}$, we get
\begin{eqnarray*}
\|\big(\int_{0}^{+\infty}\big(t^{k-\alpha_{2}}|\frac{\partial^{k}P_{t}f}{\partial
t^{k}}|\big)^{q_{2}}\frac{dt}{t}\big)^{\frac{1}{q_{2}}}\|_{p,\gamma_{d}}&\leq& C \|\big(\int_{0}^{+\infty}\big(t^{k-\alpha_{1}}|\frac{\partial^{k}P_{t}f}{\partial
t^{k}}|\big)^{q_{1}}\frac{dt}{t}\big)^{\frac{1}{q_{1}}}\|_{p,\gamma_{d}}+C \|T^{\ast}f\|_{p,\gamma_{d}}\\
&\leq&C[ \|\big(\int_{0}^{+\infty}\big(t^{k-\alpha_{1}}|\frac{\partial^{k}P_{t}f}{\partial
t^{k}}|\big)^{q_{1}}\frac{dt}{t}\big)^{\frac{1}{q_{1}}}\|_{p,\gamma_{d}} + \|f\|_{p,\gamma_{d}}]<+\infty,
\end{eqnarray*}
Thus, $f\in F_{p}^{\alpha_{2},q_{2} }(\gamma_d).$  \ep

Let us observe that the Gaussian  Besov-Lipschitz
 spaces and the Gaussian Triebel-Lizorkin spaces are by construction subspaces of $L^p(\gamma_d)$. Moreover since trivially $\left\| f \right\|_{p,\gamma_d} \leq \left\| f \right\|_{B_{p,q}^{\alpha}}$and  $\left\| f \right\|_{p,\gamma_d} \leq \left\| f \right\|_{F_{p,q}^{\alpha}}$ the inclusions are continuous.
On the other hand, from (\ref{PHHerm}) it is clear that for all $t>0$ and $k\in \mathbb{N}$, 
$$\displaystyle\frac{\partial^{k}}{\partial t^{k}}P_{t}h_{\beta}(x)=(-1)^k |\beta|^{k/2}e^{-t\sqrt{|\beta|}}h_{\beta}(x),$$
and therefore 
\begin{eqnarray*}
\big(\int_{0}^{+\infty}\big(t^{k-\alpha}\|\frac{\partial^{k}}{\partial t^{k}}P_{t}h_{\beta}\|_{p,\gamma}\big)^{q}\frac{dt}{t}\big)^{1/q}&=&\big(\int_{0}^{+\infty}
\big(t^{k-\alpha}\|(-|\beta|^{1/2})^{k} e^{-t\sqrt{|\beta|}}h_{\beta}\|_{p,\gamma}\big)^{q}\frac{dt}{t}\big)^{1/q}\\
&=& |\beta|^{k/2} \big(\int_{0}^{+\infty}
t^{(k-\alpha)q} e^{-t\sqrt{|\beta|}q}\frac{dt}{t}\big)^{1/q} \|h_{\beta}\|_{p,\gamma}\\
&=&\frac{|\beta|^{\alpha/2}}{q^{k-\alpha}} \big(\Gamma((k-\alpha)q)\big)^{1/q}  \|h_{\beta}\|_{p,\gamma}<\infty.
\end{eqnarray*}
Thus  $h_{\beta}\in B^{\alpha}_{p,q}(\gamma_{d})$ and
$$
\|h_{\beta}\|_{B^{\alpha}_{p,q}}=(1+\frac{|\beta|^{\alpha/2}}{q^{k-\alpha}}\big(\Gamma((k-\alpha)q)\big)^{1/q})\|h_{\beta}\|_{p,\gamma}.$$
Similarly,  $h_{\beta}\in F^{\alpha}_{p,q}(\gamma_{d})$ and
\begin{eqnarray*}
\|h_{\beta}\|_{F^{\alpha}_{p,q}}&=&\|h_{\beta}\|_{p,\gamma}+\|\big(\int_{0}^{+\infty}\big(t^{k-\alpha}|\frac{\partial^{k}}{\partial t^{k}}P_{t}h_{\beta}(x)|\big)^{q}\frac{dt}{t}\big)^{1/q}\|_{p,\gamma}\\
&=&(1+\frac{|\beta|^{\alpha/2}}{q^{k-\alpha}}\big(\Gamma((k-\alpha)q)\big)^{1/q})\|h_{\beta}\|_{p,\gamma} =  \|h_{\beta}\|_{B^{\alpha}_{p,q}}.
\end{eqnarray*}
Therefore, the set of polynomials ${\mathcal P}$ is included in $ B^{\alpha}_{p,q}(\gamma_{d})$ and in $ F^{\alpha}_{p,q}(\gamma_{d})$.

Also we have the following  inclusion relations between Gaussian
Triebel-Lizorkin spaces and Gaussian Besov-Lipschitz spaces,
\begin{proposition}\label{prop5}
Let $\alpha \geq 0$ and  $p,q>1$
\begin{enumerate}
\item [i)]  If $p=q$  then $$F_{p,p}^{\alpha}(\gamma_{d})=B_{p,p}^{\alpha}(\gamma_{d}).$$
\item [ii)] If $q>p$ then
$$ F_{p,q}^{\alpha}(\gamma_d)\subset B_{p,q}^{\alpha}(\gamma_d).$$
\item [iii)] If $p>q$ then
$$B_{p,q}^{\alpha}(\gamma_d) \subset F_{p,q}^{\alpha}(\gamma_d).$$
\end{enumerate}

\end{proposition}

\dem

i)  Using Tonelli's theorem, we trivially have
\begin{eqnarray*}
\|\big(\int_{0}^{+\infty}\big(t^{k-\alpha}|\frac{\partial^{k}P_t
f}{\partial
t^{k}}|\big)^{p}\frac{dt}{t}\big)^{\frac{1}{p}}\|_{p,\gamma_d}&=&\big(\int_{0}^{+\infty}t^{(k-\alpha)p}
\int_{\mathbb{R}^{d}}|\frac{\partial^{k}P_t f(x)}{\partial
t^{k}}|^{p}\gamma_{d}(dx)\frac{dt}{t}\big)^{\frac{1}{p}}\\
&=&\big(\int_{0}^{+\infty}\big(t^{k-\alpha}
\|\frac{\partial^{k}P_t f}{\partial
t^{k}}\|_{p}\big)^{p}\frac{dt}{t}\big)^{\frac{1}{p}}.
\end{eqnarray*} 

 ii) Suppose $q>p$, by Minkowski's integral
inequality we have,
\begin{eqnarray*}
\left( \int_0^{\infty} (t^{k-\alpha} \left\| \frac{\partial^{k}
P_t f}{\partial t^{k}} \right\|_{p,\gamma_d} )^{q} \frac{dt}{t}
\right) ^{p/q} &=&\left( \int_0^{\infty}
t^{(k-\alpha)q}\big(\int_{\mathbb{R}^{d}} \left|
\frac{\partial^{k}P_t f(x)}{\partial t^{k}}
\right|^{p}\gamma_{d}(dx) \big)^{q/p} \frac{dt}{t} \right) ^{p/q}\\
&\leq&\int_{\mathbb{R}^{d}}\big( \int_0^{\infty}
\big(t^{k-\alpha}\left| \frac{\partial^{k}P_t f(x)}{\partial
t^{k}}\right|\big)^{q}\frac{dt}{t} \big)^{p/q} \gamma_{d}(dx).
\end{eqnarray*}
Therefore,
\begin{eqnarray*}
\|f\|_{B_{p,q}^{\alpha}}&=&\|f\|_{p,\gamma_{d}}+\left(
\int_0^{\infty} (t^{k-\alpha} \left\| \frac{\partial^{k}P_t
f}{\partial t^{k}} \right\|_{p,\gamma_d} )^{q} \frac{dt}{t}
\right)
^{1/q}\\
&\leq&\|f\|_{p,\gamma_{d}}+\|\big( \int_0^{\infty}
\big(t^{k-\alpha}\left| \frac{\partial^{k}P_t f}{\partial
t^{k}}\right|\big)^{q}\frac{dt}{t} \big)^{1/q}
\|_{p,\gamma_{d}}=\|f\|_{F_{p,q}^{\alpha}}.
\end{eqnarray*}

 iii) Finally, if $p>q$, using  again Minkowski's integral inequality, we 
\begin{eqnarray*}
\|f\|_{F_{p,q}^{\alpha}}&=&\|f\|_{p,\gamma_{d}}+\|\big(
\int_0^{\infty} \big(t^{k-\alpha}\left| \frac{\partial^{k}P_t
f}{\partial t^{k}}\right|\big)^{q}\frac{dt}{t} \big)^{1/q}
\|_{p,\gamma_{d}}\\
&\leq&\|f\|_{p,\gamma_{d}}+\left( \int_0^{\infty} (t^{k-\alpha}
\left\| \frac{\partial^{k}P_t f}{\partial t^{k}}
\right\|_{p,\gamma_d} )^{q} \frac{dt}{t} \right)
^{1/q}=\|f\|_{B_{p,q}^{\alpha}}.
\end{eqnarray*}    \ep\\

Let us prove now that the Gaussian Sobolev spaces
$L^p_\alpha(\gamma_d)$ are contained in some Besov-Lipschitz and
Triebel-Lizorkin spaces, and therefore they are ``finer scales" to
measure the regularity of functions.

\begin{theorem}
Let us suppose that $1<p<+\infty$ and $\alpha > 0$. Then
\begin{enumerate}
\item [i)] $ L^p_\alpha(\gamma_d)\subset F_{p,2}^{\alpha}(\gamma_{d})$ if $p>1$.
\item [ii)] $ L^p_\alpha(\gamma_d)\subset
B_{p,p}^{\alpha}(\gamma_{d})=F_{p,p}^{\alpha}(\gamma_{d})$ if $p\geq 2$. 
\item [iii)] $
L^p_\alpha(\gamma_d)\subset B_{p,2}^{\alpha}(\gamma_{d})$ if $p\leq 2$. 

\end{enumerate}
\end{theorem}

\dem
 
 i). We have to consider two cases:
\begin{enumerate}
\item If $\alpha\geq 1$.
Suppose $h\in L^p_\alpha(\gamma_d)$
then $h=\mathcal{J}_{\alpha}f$, $f\in L^{p}(\gamma_{d})$, by the
change of variable $u=t+s$ using  the fact the representation of the Bessel potentials (\ref{Beselrepre}) and  Hardy's
inequality to get,

\begin{eqnarray*}
\big(\int_{0}^{+\infty}\big(t^{k-\alpha}|\frac{\partial^{k}P_{t}h(x)}{\partial
t^{k}}|\big)^{2}\frac{dt}{t}\big)^{\frac{1}{2}}&=&\big(\int_{0}^{+\infty}t^{2(k-\alpha)}|\frac{\partial^{k}P_{t}\mathcal{J}_{\alpha}f(x)}{\partial
t^{k}}|^{2}\frac{dt}{t}\big)^{\frac{1}{2}}\\
&\leq&\frac{1}{\Gamma(\alpha)}\big(\int_{0}^{+\infty}t^{2(k-\alpha)}\big(\int_{0}^{+\infty}s^{\alpha}e^{-s}
|\frac{\partial^{k}P_{t+s}f(x)}{\partial(t+s)^{k}}|\frac{ds}{s}\big)^{2}\frac{dt}{t}\big)^{\frac{1}{2}}\\
&=&\frac{1}{\Gamma(\alpha)}\big(\int_{0}^{+\infty}t^{2(k-\alpha)}\big(\int_{t}^{+\infty}(u-t)^{\alpha-1}e^{t-u}
|\frac{\partial^{k}P_{u}f(x)}{\partial u^{k}}|du\big)^{2}\frac{dt}{t}\big)^{\frac{1}{2}}\\
&\leq&\frac{1}{\Gamma(\alpha)}\big(\int_{0}^{+\infty}\big(\int_{t}^{+\infty}u^{\alpha-1}
|\frac{\partial^{k}P_{u}f(x)}{\partial
u^{k}}|du\big)^{2}t^{2(k-\alpha)-1}dt\big)^{\frac{1}{2}}\\
&\leq& \frac{1}{\Gamma(\alpha)}\frac{1}{k-\alpha}
\big(\int_{0}^{+\infty}\big(u^{k}
|\frac{\partial^{k}P_{u}f(x)}{\partial
u^{k}}|\big)^{2}\frac{du}{u}\big)^{\frac{1}{2}}.
\end{eqnarray*}

Hence,  by the $L^p(\gamma_d)$-continuity of the Gaussian
Littlewood-Paley  $g_k$- function, see \cite{forscotur}
\begin{eqnarray*}
\|\big(\int_{0}^{+\infty}\big(t^{k-\alpha}|\frac{\partial^{k}P_{t}h}{\partial
t^{k}}|\big)^{2}\frac{dt}{t}\big)^{\frac{1}{2}}\|_{p,\gamma}&\leq&\frac{1}{\Gamma(\alpha)}\frac{1}{k-\alpha}
\|\big(\int_{0}^{+\infty}\big(u^{k}
|\frac{\partial^{k}P_{u}f}{\partial
u^{k}}|\big)^{2}\frac{du}{u}\big)^{\frac{1}{2}}\|_{p,\gamma}\\
&=&C_{k,\alpha} \|g_{k}f\|_{p,\gamma} \leq C_{k,\alpha}
\|f\|_{p,\gamma}=C_{k,\alpha} \|h\|_{p,\alpha},
\end{eqnarray*}
thus $h\in  F_{p,2}^{\alpha}(\gamma_d)$.\\

\item If $0\leq \alpha< 1$.
 Suppose
$h\in L^p_\alpha(\gamma_d)$, then
$h={\mathcal{J}}_{\alpha}f$, $f\in L^{p}(\gamma_{d}),$ again using (\ref{Beselrepre}),
\begin{eqnarray*}
\big(\int_{0}^{+\infty}\big(t^{k-\alpha}|\frac{\partial^{k}P_{t}h(x)}{\partial
t^{k}}|\big)^{2}\frac{dt}{t}\big)^{\frac{1}{2}}
&\leq&\frac{1}{\Gamma(\alpha)}\big(\int_{0}^{+\infty}t^{2(k-\alpha)}\big(\int_{0}^{+\infty}s^{\alpha}e^{-s}
|\frac{\partial^{k}P_{t+s}f(x)}{\partial(t+s)^{k}}|\frac{ds}{s}\big)^{2}\frac{dt}{t}\big)^{\frac{1}{2}} \\
&\leq&\frac{C}{\Gamma(\alpha)}\big(\int_{0}^{+\infty}t^{2(k-\alpha)-1}\big[(\int_{0}^{t}s^{\alpha}e^{-s}
|\frac{\partial^{k}P_{t+s}f(x)}{\partial(t+s)^{k}}|\frac{ds}{s})^{2}\\
&& \quad \quad \quad  \quad \quad \quad +(\int_{t}^{+\infty}s^{\alpha}e^{-s}
|\frac{\partial^{k}P_{t+s}f(x)}{\partial(t+s)^{k}}|\frac{ds}{s})^{2}\big]dt\big)^{\frac{1}{2}}\\
&\leq&\frac{C}{\Gamma(\alpha)}(\int_{0}^{+\infty}t^{2(k-\alpha)-1}\big(\int_{0}^{t}s^{\alpha-1}e^{-s}
|\frac{\partial^{k}P_{t+s}f(x)}{\partial(t+s)^{k}}|ds\big)^{2}dt\big)^{\frac{1}{2}}\\
&& \quad +\frac{C}{\Gamma(\alpha)} (\int_{0}^{+\infty}t^{2(k-\alpha)-1}\big(\int_{t}^{+\infty}s^{\alpha-1}e^{-s}
|\frac{\partial^{k}P_{t+s}f(x)}{\partial(t+s)^{k}}|ds\big)^{2}dt\big)^{\frac{1}{2}}\\
&=& I + II.
\end{eqnarray*}
Now, since $e^{-s}<1$,  $s^{\alpha-1}<t^{\alpha-1}$ as $\alpha<1$, and
using the change of variables $u=t+s$ and Hardy inequality we get,
\begin{eqnarray*}
  II &\leq&\big(\int_{0}^{+\infty}t^{2(k-1)-1}\big(\int_{t}^{+\infty}
|\frac{\partial^{k}P_{t+s}f(x)}{\partial(t+s)^{k}}|ds\big)^{2}dt\big)^{\frac{1}{2}} =\big(\int_{0}^{+\infty}t^{2(k-1)-1}\big(\int_{2t}^{+\infty}
|\frac{\partial^{k}P_{u}f(x)}{\partial u^{k}}|du\big)^{2}dt\big)^{\frac{1}{2}}\\
&\leq&\big(\int_{0}^{+\infty}t^{2(k-1)-1}\big(\int_{t}^{+\infty}
|\frac{\partial^{k}P_{u}f(x)}{\partial u^{k}}|du\big)^{2}dt\big)^{\frac{1}{2}}
\leq \big(\int_{0}^{+\infty}\big(u |\frac{\partial^{k}P_{u}f(x)}{\partial u^{k}}|\big)^{2}u^{2(k-1)-1}du\big)^{\frac{1}{2}}.\\
&=&\big(\int_{0}^{+\infty}\big(|u^{k}\frac{\partial^{k}P_{u}f(x)}{\partial u^{k}}|\big)^{2}\frac{du}{u}\big)^{\frac{1}{2}}=g_{k}f(x).
\end{eqnarray*}
On the other hand, again since $e^{-s}<1$,
\begin{eqnarray*}
I^2 &\leq&\int_{0}^{+\infty}t^{2(k-\alpha)-1}\big(\int_{0}^{t}s^{\alpha-1}
|\frac{\partial^{k}P_{t+s}f(x)}{\partial(t+s)^{k}}|ds\big)^{2}dt=\frac{1}{\alpha^{2}}\int_{0}^{+\infty}t^{2k-1}\big(\frac{\alpha
}{ t^{\alpha}}\int_{0}^{t}s^{\alpha-1}
|\frac{\partial^{k}P_{t+s}f(x)}{\partial(t+s)^{k}}|ds\big)^{2}dt
\end{eqnarray*}
Now, as $\alpha> 0$ using Jensen's inequality for the measure $\frac{\alpha}{t^{\alpha}}s^{\alpha-1} ds$ and
Tonelli's Theorem,
\begin{eqnarray*}
I^2&\leq&\frac{1}{\alpha^{2}}\int_{0}^{+\infty}t^{2k-1}\big(\frac{\alpha
}{ t^{\alpha}}\int_{0}^{t}s^{\alpha-1}
|\frac{\partial^{k}P_{t+s}f(x)}{\partial(t+s)^{k}}|^{2}ds\big)dt\\
&\leq&\frac{1}{\alpha}\int_{0}^{+\infty}s^{\alpha-1}\big(\int_{s}^{+\infty}(t+s)^{2k-\alpha-1}
|\frac{\partial^{k}P_{t+s}f(x)}{\partial(t+s)^{k}}|^{2}dt\big)ds,
\end{eqnarray*}
since $2k-\alpha-1>0.$ Finally, again using the change of variables $u=t+s$ and the Hardy inequality
\begin{eqnarray*}
I^2&\leq&\frac{1}{\alpha}\int_{0}^{+\infty}s^{\alpha-1}\big(\int_{2s}^{+\infty}u^{2k-\alpha-1}
|\frac{\partial^{k}P_{u}f(x)}{\partial u^{k}}|^{2}du\big)ds\\
&\leq& \frac{1}{\alpha}\int_{0}^{+\infty}s^{\alpha-1}\big(\int_{s}^{+\infty}u^{2k-\alpha-1}
|\frac{\partial^{k}P_{u}f(x)}{\partial u^{k}}|^{2}du\big)ds\\
&\leq&\frac{1}{\alpha}\int_{0}^{+\infty}\big(u^{k}|\frac{\partial^{k}P_{u}f(x)}{\partial
u^{k}}|\big)^{2}\frac{du}{u} = \frac{1}{\alpha} g^2_{k}f(x).
\end{eqnarray*}
Hence, again by the $L^p(\gamma_d)$-continuity of the Gaussian
Littlewood-Paley  $g_k$- function,
\begin{eqnarray*}
\|\big(\int_{0}^{+\infty}\big(t^{k-\alpha}|\frac{\partial^{k}P_{t}h}{\partial
t^{k}}|\big)^{2}\frac{dt}{t}\big)^{\frac{1}{2}}\|_{p,\gamma}
\leq C_{k,\alpha}\|g_{k}f\|_{p,\gamma}\leq C_{k,\alpha}\|f\|_{p,\gamma}
=C_{k,\alpha}\|h\|_{p,\alpha}.
\end{eqnarray*}

Thus $h\in  F_{p,2}^{\alpha}(\gamma_d)$, for $0<\alpha<1$.
\end{enumerate}

Let us prove now ii). Suppose $h\in L^p_\alpha(\gamma_d)$ with
$p\geq 2$ then $h=\mathcal{J}_{\alpha}f$, $f\in
L^{p}(\gamma_{d})$. Using the inequality $(a+b)^{p}\leq
C_{p}(a^{p}+b^{p})$ if $a,b\geq 0, p\geq 1$ we get

\begin{eqnarray*}
&&\big(\int_{0}^{+\infty}\big(t^{k-\alpha}\|\frac{\partial^{k}
P_{t}\mathcal{J}_{\alpha}f}{\partial
t^{k}}\|_{p,\gamma_d}\big)^{p}\frac{dt}{t}\big)^{\frac{1}{p}}
\quad \quad\quad \quad \quad \quad \quad \quad\quad\quad
\quad\quad \quad \quad \quad \\
&&\quad \quad\quad \quad \quad\leq \frac{1}{\Gamma(\alpha)}
\big(\int_{0}^{+\infty}\big(t^{k-\alpha}\int_{0}^{+\infty}s^{\alpha}e^{-s}
\|\frac{\partial^{k}P_{t+s}f}{\partial(t+s)^{k}}\|_{p,\gamma}\frac{ds}{s}\big)^{p}\frac{dt}{t}\big)^{\frac{1}{p}}\\
&&\quad \quad\quad \quad
\quad\leq\frac{C}{\Gamma(\alpha)}\big(\int_{0}^{+\infty}t^{p(k-\alpha)}\big(\int_{0}^{t}s^{\alpha}
\|\frac{\partial^{k}P_{s+t}f}{\partial(s+t)^{k}}\|_{p,\gamma}\frac{ds}{s}\big)^{p}\\
&&\quad \quad \quad \quad \quad  \quad \quad\quad \quad  \quad
\quad\quad \quad \quad + \big(\int_{t}^{+\infty}s^{\alpha}
\|\frac{\partial^{k}P_{s+t}f}{\partial(s+t)^{k}}\|_{p,\gamma}\frac{ds}{s}\big)^{p}\frac{dt}{t}\big)^{\frac{1}{p}}.
\end{eqnarray*}

Using the inequality $(a+b)^{1/p}\leq a^{1/p}+b^{1/p}$ if $a,b\geq 0,
p\geq 1$
\begin{eqnarray*}
&&
\frac{C}{\Gamma(\alpha)}\big(\int_{0}^{+\infty}t^{p(k-\alpha)}\big(\int_{0}^{t}s^{\alpha}
\|\frac{\partial^{k}P_{s+t}f}{\partial(s+t)^{k}}\|_{p,\gamma}\frac{ds}{s}\big)^{p}\\
&&\quad \quad \quad \quad \quad  \quad \quad\quad \quad  \quad
\quad\quad \quad \quad + \big(\int_{t}^{+\infty}s^{\alpha}
\|\frac{\partial^{k}P_{s+t}f}{\partial(s+t)^{k}}\|_{p,\gamma}\frac{ds}{s}\big)^{p}\frac{dt}{t}\big)^{\frac{1}{p}}\\
&\leq&\frac{C}{\Gamma(\alpha)}\big(\int_{0}^{+\infty}t^{(k-\alpha)p}\big(\int_{0}^{t}s^{\alpha}
\|\frac{\partial^{k}P_{s+t}f}{\partial(s+t)^{k}}\|_{p,\gamma}\frac{ds}{s}\big)^{p}\frac{dt}{t}\big)^{\frac{1}{p}}\\
&& \quad \quad \quad  \quad \quad\quad
+\frac{C}{\Gamma(\alpha)}\big(\int_{0}^{+\infty}t^{(k-\alpha)p}\big(\int_{t}^{+\infty}s^{\alpha}
\|\frac{\partial^{k}P_{s+t}f}{\partial(s+t)^{k}}\|_{p,\gamma}\frac{ds}{s}\big)^{p}\frac{dt}{t}\big)^{\frac{1}{p}}\\
&=&I+II.
\end{eqnarray*}

Now,  using again the Hardy's  inequality, since $k>\alpha$ and
lemma \ref{kdecay}
\begin{eqnarray*}
II&=&\frac{C}{\Gamma(\alpha)}\big(\int_{0}^{+\infty}t^{p(k-\alpha)}\big(\int_{t}^{+\infty}s^{\alpha}
\|\frac{\partial^{k}P_{s+t}f}{\partial(s+t)^{k}}\|_{p,\gamma}\frac{ds}{s}\big)^{p}\frac{dt}{t}\big)^{\frac{1}{p}}\\
&\leq&\frac{C}{\Gamma(\alpha)}\big(\int_{0}^{+\infty}t^{p(k-\alpha)}\big(\int_{t}^{+\infty}s^{\alpha}
\|\frac{\partial^{k}P_{s}f}{\partial s^{k}}\|_{p,\gamma}\frac{ds}{s}\big)^{p}\frac{dt}{t}\big)^{\frac{1}{p}}\\
&\leq&\frac{C}{\Gamma(\alpha)}\frac{1}{k-\alpha}\big(\int_{0}^{+\infty}\big(s^{\alpha}
\|\frac{\partial^{k}}{\partial s^{k}}P_{s}f\|_{p,\gamma}\big)^{p}s^{(k-\alpha)p-1}ds\big)^{\frac{1}{p}}\\
&=&C_{k,\alpha}\big(\int_{0}^{+\infty}\big(s^{k}
\|\frac{\partial^{k}}{\partial
s^{k}}P_{s}f\|_{p,\gamma}\big)^{p}\frac{ds}{s}\big)^{\frac{1}{p}}=C_{k,\alpha}\|\big(
\int_{0}^{+\infty}|s^{k}\frac{\partial^{k}P_{s}f}{\partial
s^{k}}|^{p}\frac{ds}{s}\big)^{\frac{1}{p}}\|_{p},
\end{eqnarray*}
by Tonelli's Theorem.

Now since $p\geq 2$ using Lemma \ref{lematec5}, we have

\begin{eqnarray*}
\int_{0}^{+\infty}|u^{k}\frac{\partial^{k}P_{u}f(x)}{\partial
u^{k}}|^{p}\frac{du}{u}&=&\int_{0}^{+\infty}\big(u^{k}|\frac{\partial^{k}}{\partial
u^{k}}P_{u}f(x)|\big)^{p-2}
\big(u^{k}|\frac{\partial^{k}}{\partial
u^{k}}P_{u}f(x)|\big)^{2}\frac{du}{u}\\
&\leq&C\big(T^{\ast}f(x)\big)^{p-2}
\int_{0}^{+\infty}\big(u^{k}|\frac{\partial^{k}}{\partial
u^{k}}P_{u}f(x)|\big)^{2}\frac{du}{u}.
\end{eqnarray*}
Therefore

\begin{eqnarray*}
\|\big( \int_{0}^{+\infty}|u^{k}\frac{\partial^{k}P_{u}f}{\partial
u^{k}}|^{p}\frac{du}{u}\big)^{\frac{1}{p}}\|_{p}^{p}&=&\int_{\mathbb{R}^{d}}\big(
\int_{0}^{+\infty}|u^{k}\frac{\partial^{k}P_{u}f(x)}{\partial
u^{k}}|^{p}\frac{du}{u}\big)\gamma_{d}(dx)\\
&\leq&C\int_{\mathbb{R}^{d}}\big(\big(T^{\ast}f(x)\big)^{p-2}
\int_{0}^{+\infty}\big(u^{k}|\frac{\partial^{k}P_{u}f(x)}{\partial
u^{k}}|\big)^{2}\frac{du}{u})\gamma_{d}(dx)
\end{eqnarray*}

Using  H\"older inequality, with
$\theta=\displaystyle\frac{2}{p}$,  and the $L^{p}(\gamma_{d})$
continuity of $T^{\ast}$ and $g_{k}$, we have
\begin{eqnarray*}
\|\big( \int_{0}^{+\infty}|u^{k}\frac{\partial^{k}P_{u}f}{\partial
u^{k}}|^{p}\frac{du}{u}\big)^{\frac{1}{p}}\|_{p}^{p}&\leq&C\int_{\mathbb{R}^{d}}\big(\big(T^{\ast}f(x)\big)^{p-2}
\int_{0}^{+\infty}\big(u^{k}|\frac{\partial^{k}P_{u}f(x)}{\partial
u^{k}}|\big)^{2}\frac{du}{u})\gamma_{d}(dx)\\
&\leq&C\big(\int_{\mathbb{R}^{d}}\big(\big(T^{\ast}f(x)\big)^{(p-2).\frac{1}{1-\theta}}\gamma_{d}(dx)\big)^{1-\theta}.\\
&&\quad \quad
\times\big(\int_{\mathbb{R}^{d}}\big(\int_{0}^{+\infty}\big(u^{k}|\frac{\partial^{k}}{\partial
u^{k}}P_{u}f(x)|\big)^{2}\frac{du}{u}\big)^{\frac{1}{\theta}}\gamma_{d}(dx)\big)^{\theta}\\
&=&C\big(\int_{\mathbb{R}^{d}}\big(\big(T^{\ast}f(x)\big)^{p}\gamma_{d}(dx)\big)^{\frac{p-2}{p}}.\\
&&\quad \quad \times
\big(\int_{\mathbb{R}^{d}}\big(\int_{0}^{+\infty}\big(u^{k}|\frac{\partial^{k}}{\partial
u^{k}}P_{u}f(x)|\big)^{2}\frac{du}{u}\big)^{\frac{p}{2}}\gamma_{d}(dx)\big)^{\frac{2}{p}}\\
&=&C\|T^{\ast}f\|_{p,\gamma}^{p-2} \|g_{k}f\|_{p,\gamma}^{2} \leq
C\|f\|_{p,\gamma}^{p}.
\end{eqnarray*}
Thus,
$$
II\leq C_{k,\alpha} \|h\|_{p,\alpha}.
$$
Now, using again lemma \ref{kdecay}  and since $\alpha>0$
\begin{eqnarray*}
I&=&\frac{C}{\Gamma(\alpha)}\big(\int_{0}^{+\infty}t^{p(k-\alpha)}\big(\int_{0}^{t}s^{\alpha}
\|\frac{\partial^{k}}{\partial(s+t)^{k}}P_{s+t}f\|_{p,\gamma}\frac{ds}{s}\big)^{p}\frac{dt}{t}\big)^{\frac{1}{p}}\\
&\leq&\frac{C}{\Gamma(\alpha)}\big(\int_{0}^{+\infty}t^{p(k-\alpha)}\big(\int_{0}^{t}s^{\alpha}
\|\frac{\partial^{k}P_{t}f}{\partial t^{k}}\|_{p,\gamma}\frac{ds}{s}\big)^{p}\frac{dt}{t}\big)^{\frac{1}{p}}\\
&=&\frac{1}{\alpha}\frac{C}{\Gamma(\alpha)}\big(\int_{0}^{+\infty}t^{k}\|\frac{\partial^{k}P_{t}f}{\partial t^{k}}\|_{p,\gamma}^{p}\frac{dt}{t}\big)^{\frac{1}{p}}
\leq C_{k,\alpha} \|h\|_{p,\alpha},
\end{eqnarray*}
So $h\in  B_{p,p}^{\alpha}(\gamma_d)$, if $p\geq 2$.\\

 iii) can be proved using similar arguments as in i) and ii) but it is  immediate consequences of  i) and 
 of Proposition \ref{prop5} ii). \ep\\

In \cite{iris}, using Theorem 3.2,  it is claimed that the  Gaussian Sobolev spaces $L_{\alpha}^p(\gamma_d)$ coincide with the homogeneous Gaussian Triebel-Lizorkin $\dot{F}_{p,2}^{\alpha}$ but the proof of that theorem  is wrong since it is assumed the the operator involved is linear but it is actually only sublinear.\\

Now, let us prove some interpolation results for the Gaussian
Besov-Lipschitz spaces and for the Gaussian Triebel-Lizorkin
spaces,
\begin{theorem} We have the following interpolation results:
\begin{enumerate}
\item [i)] For $1<p_{j},q_{j}<+\infty$ and $\alpha_{j}\geq0$, if
$f\in B_{p_{j},q_{j}}^{\alpha_{j}}(\gamma_d)$, $j=0,1,$ then $f\in
B_{p,q}^{\alpha}(\gamma_d)$, where
$\alpha=\alpha_{0}(1-\theta)+\alpha_{1}\theta$, and
$$\displaystyle\frac{1}{p}=\frac{1}{p_{0}}(1-\theta)+\frac{\theta}{p_{1}}, \frac{1}{q}=\frac{1}{q_{0}}(1-\theta)+\frac{\theta}{q_{1}}, \, \, 0<\theta<1.$$

\item[ii)] For $1<p_{j},q_{j}<+\infty$ and $\alpha_{j}\geq0$, if
$f\in F_{p_{j},q_{j}}^{\alpha_{j}}(\gamma_d)$, $j=0,1,$ then $f\in
F_{p,q}^{\alpha}(\gamma_d)$, where
$\alpha=\alpha_{0}(1-\theta)+\alpha_{1}\theta$, and
$$\displaystyle\frac{1}{p}=\frac{1}{p_{0}}(1-\theta)+\frac{\theta}{p_{1}}, \frac{1}{q}=\frac{1}{q_{0}}(1-\theta)+\frac{\theta}{q_{1}}, \, \, 0<\theta<1.$$
\end{enumerate}
\end{theorem}
\dem The proof of both results are based  in the following
interpolation result  for $L^p(\gamma_d)$ spaces (actually true
for any measure $\mu$) that is obtained using H\"older inequality:

For $1<r_{0},r_{1}<\infty$ and
$\displaystyle\frac{1}{r}=\frac{1}{r_{0}}(1-\eta)+\frac{\eta}{r_{1}},
0<\eta<1$. If $f\in L^{r_{j}}(\gamma_d)$, $j=0,1$ then $f\in
L^{r}(\gamma_d)$ and
\begin{equation}
\|f\|_{r,\gamma_d}\leq\|f\|_{r_{0},\gamma_d}^{1-\eta}\|f\|_{r_{1},\gamma_d}^{\eta}.
\end{equation}

Let us prove i).  Let $k$ be any integer greater than
$\alpha_{0}$ and $\alpha_{1}$, by using the above result we get
for $\alpha = \alpha_0 (1- \theta) +\alpha_1 \theta$,

\begin{eqnarray*}
\int_{0}^{+\infty}\big(t^{k-\alpha}\|\frac{\partial^{k}P_{t}f}{\partial
t^{k}}\|_{p,\gamma_d}\big)^{q}\frac{dt}{t} &\leq& \int_{0}^{+\infty}\big(t^{k-(\alpha_0 (1- \theta) +\alpha_1
\theta)}\|\frac{\partial^{k}P_{t}f}{\partial
t^{k}}\|_{p_{0},\gamma_d}^{1-\theta}\|\frac{\partial^{k}P_{t}f}{\partial
t^{k}}\|_{p_{1},\gamma_d}^{\theta}\big)^{q}\frac{dt}{t}\\
&=&\int_{0}^{+\infty}\big(t^{(1-\theta)(k-\alpha_{0})+\theta(k-\alpha_{1})}\|\frac{\partial^{k}P_{t}f}{\partial
t^{k}}\|_{p_{0},\gamma_d}^{1-\theta}\|\frac{\partial^{k}P_{t}f}{\partial
t^{k}}\|_{p_{1},\gamma_d}^{\theta}\big)^{q}\frac{dt}{t}\\
&=&\int_{0}^{+\infty}\big(t^{k-\alpha_{0}}\|\frac{\partial^{k}P_{t}f}{\partial
t^{k}}\|_{p_{0},\gamma_d}\big)^{(1-\theta)q}\big(t^{k-\alpha_{1}}\|\frac{\partial^{k}P_{t}f}{\partial
t^{k}}\|_{p_{1},\gamma_d}\big)^{\theta q}\frac{dt}{t}.
\end{eqnarray*}
Now, if $\displaystyle\lambda=\frac{\theta q}{q_{1}}$ then
$0<\lambda<1$ and $q=(1-\lambda)q_{0}+\lambda q_{1}$. Therefore by
using again the H\"older inequality ,
\begin{eqnarray*}
&&\int_{0}^{+\infty}\big(t^{k-\alpha}\|\frac{\partial^{k}P_{t}f}{\partial
t^{k} }\|_{p,\gamma_d}\big)^{q}\frac{dt}{t} \quad \quad\quad \quad
\quad \quad \quad \quad\quad \quad \quad\quad \quad\quad \quad
\quad\\
&&\quad \quad\quad \quad \leq \big(\int_{0}^{+\infty}\big(t^{k-\alpha_{0}}\|\frac{\partial^{k}P_{t}f}{\partial
t^{k}}\|_{p_{0},\gamma_d}\big)^{
q_{0}}\frac{dt}{t}\big)^{1-\lambda}
\big(\int_{0}^{+\infty}\big(t^{k-\alpha_{1}}\|\frac{\partial^{k}P_{t}f}{\partial
t^{k}}\|_{p_{1},\gamma_d}\big)^{
q_{1}}\frac{dt}{t}\big)^{\lambda}<\infty,
\end{eqnarray*}
and so $f\in B_{p,q}^{\alpha}(\gamma_d)$.

ii) Analogously, by taking $\beta=\frac{p\theta}{p_{1}}, \, \lambda=\frac{q\theta}{q_{1}},
$ we have $0<\beta,\lambda<1$ and
$p=(1-\beta)p_{0}+\beta p_{1}, q=(1-\lambda)q_{0}+\lambda q_{1}$.
Let $k$ be any integer greater than $\alpha_{0}$ and $\alpha_{1}$,
by using H\"older we get for $\alpha = \alpha_0 (1- \theta)
+\alpha_1 \theta$,

\begin{eqnarray*}
\int_{0}^{+\infty}\big(t^{k-\alpha}|\frac{\partial^{k}P_{t}f}{\partial
t^{k} }|\big)^{q}\frac{dt}{t} &=&\int_{0}^{+\infty}\big(t^{k-\alpha_{0}}
|\frac{\partial^{k}P_{t}f}{\partial t^{k}
}|\big)^{(1-\theta)q}\big(t^{k-\alpha_{1}}|\frac{\partial^{k}P_{t}f}{\partial
t^{k}}|\big)^{\theta q}\frac{dt}{t}\\
&=&\int_{0}^{+\infty}\big(t^{k-\alpha_{0}}
|\frac{\partial^{k}P_{t}f}{\partial t^{k}
}|\big)^{(1-\lambda)q_{0}}\big(t^{k-\alpha_{1}}|\frac{\partial^{k}P_{t}f}{\partial
t^{k}}|\big)^{\lambda q_{1}}\frac{dt}{t}\\
&\leq&\big(\int_{0}^{+\infty}\big(t^{k-\alpha_{0}}
|\frac{\partial^{k}P_{t}f}{\partial
t^{k}}|\big)^{q_{0}}\frac{dt}{t}\big)^{1-\lambda}\big(\int_{0}^{+\infty}\big(t^{k-\alpha_{1}}|\frac{\partial^{k}P_{t}f}{\partial
t^{k} }| \big)^{q_{1}}\frac{dt}{t}\big)^{\lambda}.
\end{eqnarray*}

Thus,
\begin{eqnarray*}
&&\|\big(\int_{0}^{+\infty}\big(t^{k-\alpha}|\frac{\partial^{k}P_{t}f}{\partial
t^{k}}|\big)^{q}\frac{dt}{t}\big)^{\frac{1}{q}}\|_{p,\gamma_{d}}^{p}
=   \int_{\mathbb{R}^{d}}\big(\int_{0}^{+\infty}\big(t^{k-\alpha}
|\frac{\partial^{k}P_{t}f}{\partial
t^{k}}|\big)^{q}\frac{dt}{t}\big)^{\frac{p}{q}}\gamma_{d}(dx)\\
&& \quad \quad \quad \leq
\int_{\mathbb{R}^{d}}\big(\int_{0}^{+\infty}\big(t^{k-\alpha_{0}}
|\frac{\partial^{k}P_{t}f}{\partial
t^{k}}|\big)^{q_{0}}\frac{dt}{t}\big)^{\frac{(1-\lambda)p}{q}}
\big(\int_{0}^{+\infty}\big(t^{k-\alpha_{1}}|\frac{\partial^{k}P_{t}f}{\partial
t^{k} }|\big)^{q_{1}}\frac{dt}{t}\big)^{\frac{\lambda p}{q}}\gamma_{d}(dx)\\
&& \quad \quad \quad =\int_{\mathbb{R}^{d}}\big(\int_{0}^{+\infty}\big(t^{k-\alpha_{0}}
|\frac{\partial^{k}P_{t}f}{\partial
t^{k}}|\big)^{q_{0}}\frac{dt}{t}\big)^{\frac{(1-\theta)p}{q_{0}}}
\big(\int_{0}^{+\infty}\big(t^{k-\alpha_{1}}|\frac{\partial^{k}P_{t}f}{\partial t^{k} }|\big)^{q_{1}}\frac{dt}{t}\big)^{\frac{\theta p}{q_{1}}}\gamma_{d}(dx)\\
&& \quad \quad \quad =\int_{\mathbb{R}^{d}}\big(\int_{0}^{+\infty}\big(t^{k-\alpha_{0}}
|\frac{\partial^{k}P_{t}f}{\partial
t^{k}}|\big)^{q_{0}}\frac{dt}{t}\big)^{\frac{(1-\beta)p_{0}}{q_{0}}}\big(\int_{0}^{+\infty}\big(t^{k-\alpha_{1}}|\frac{\partial^{k}P_{t}f}{\partial
t^{k} }|\big)^{q_{1}}\frac{dt}{t}\big)^{\frac{\beta
p_{1}}{q_{1}}}\gamma_{d}(dx),
\end{eqnarray*}
and then again using H\"older inequality,
\begin{eqnarray*}
&&\|\big(\int_{0}^{+\infty}\big(t^{k-\alpha}|\frac{\partial^{k}
P_{t}f}{\partial
t^{k}}|\big)^{q}\frac{dt}{t}\big)^{\frac{1}{q}}\|_{p,\gamma_{d}}^{p}
\quad \quad\quad \quad \quad \quad \quad \quad\quad \quad \quad
\quad \\
&&\quad\quad \quad \quad \leq
\big(\int_{\mathbb{R}^{d}}\big(\int_{0}^{+\infty}\big(t^{k-\alpha_{0}}
|\frac{\partial^{k}P_{t}f}{\partial t^{k}}|\big)^{q_{0}}\frac{dt}{t}\big)^{\frac{p_{0}}{q_{0}}}\gamma_{d}(dx)\big)^{1-\beta}\\
&&\quad \quad\quad \quad \quad \quad \quad\quad \quad\quad \quad
\times\big(\int_{\mathbb{R}^{d}}\big(\int_{0}^{+\infty}\big(t^{k-\alpha_{1}}|\frac{\partial^{k}P_{t}f}{\partial
t^{k} }|\big)^{q_{1}}\frac{dt}{t}\big)^{\frac{
p_{1}}{q_{1}}}\gamma_{d}(dx)\big)^{\beta}\\
&=&\|\big(\int_{0}^{+\infty}\big(t^{k-\alpha_{0}}
|\frac{\partial^{k}P_{t}f}{\partial t^{k}}|\big)^{q_{0}}\frac{dt}{t}\big)^{\frac{1}{q_{0}}}\|_{p_{0},\gamma_{d}}^{p_{0}(1-\beta)}\\
&&\quad \quad\quad \quad \quad \quad \quad\quad \quad\quad \quad
\times\|\big(\int_{0}^{+\infty}\big(t^{k-\alpha_{1}}
|\frac{\partial^{k}P_{t}f}{\partial t^{k}
}|\big)^{q_{1}}\frac{dt}{t}\big)^{\frac{1}{q_{1}}}\|_{p_{1},\gamma_{d}}^{p_{1}\beta}<
+\infty.
\end{eqnarray*}
Hence $f\in F_{p,q}^{\alpha}(\gamma_d)$. \ep\\

Now, we are going to study the continuity properties of the
Ornstein-Uhlenbeck semigroup,
 the Poisson-Hermite semigroup and the Bessel potentials on the Besov-Lipschitz and Triebel-Lizorkin spaces,

\begin{theorem} For the  Besov-Lipschitz spaces $B_{p,q}^{\alpha}(\gamma_{d})$ and Triebel-Lizorkin spaces $F_{p,q}^{\alpha}(\gamma_{d})$, we have
\begin{enumerate}
\item[i)]  The Ornstein-Uhlenbeck semigroup $\{T_{t}\}$ and the
Poisson-Hermite semigroup $\{P_{t}\}$ are bounded on
$B_{p,q}^{\alpha}(\gamma_{d})$.
\item[ii)] For every $\beta >0$,
the Bessel potentials $\mathcal{J}_\beta^{\gamma}$ are bounded  on
$B_{p,q}^{\alpha}(\gamma_{d})$.
\item[iii)] The Ornstein-Uhlenbeck
semigroup $\{T_{t}\}$, the Poisson-Hermite semigroup
 $\{P_{t}\}$ are bounded on $F_{p,q}^{\alpha}$.
 \item[iv)] the Bessel potentials  $\mathcal{J}_\beta^{\gamma}$ are bounded on $F_{p,q}^{\alpha}(\gamma_{d})$.
 \end{enumerate}
 \end{theorem}

\dem
\begin{enumerate}
\item[i)]  Let us prove the $B_{p,q}^{\alpha}(\gamma_{d})$-continuity of
$P_t$ for any $t >0$, the proof for $T_t$ is totally analogous. By
the $L^p$-continuity of the Poisson-Hermite semigroup, the Lebesgue's
dominated convergence theorem and Jensen's inequality we get
\begin{eqnarray*}
\int_{\mathbb{R}^{d}}| \frac{\partial^{k}
P_{t}\big(P_{s}f\big)}{\partial t^{k}}(x)|^{p}\gamma_{d}(dx)&=&
\int_{\mathbb{R}^{d}}|P_{s}
\big(\frac{\partial^{k}P_{t}f}{\partial
t^{k}}\big)(x)|^{p}\gamma_{d}(dx)\\
&\leq&\int_{\mathbb{R}^{d}}P_{s}(|\frac{\partial^{k}P_{t}f(x)}{\partial
t^{k}}|^{p})\gamma_{d}(dx)=\int_{\mathbb{R}^{d}}|\frac{\partial^{k}P_{t}f(x)}{\partial
t^{k}}|^{p}\gamma_{d}(dx).
\end{eqnarray*}
Thus,
$$\displaystyle\|\frac{\partial^{k}P_{t}(P_{s}f)}{\partial
t^{k}}\|_{p,\gamma_{d}}\leq\|\frac{\partial^{k}P_{t}f}{\partial
t^{k}}\|_{p,\gamma_{d}},$$
 and therefore
\begin{eqnarray*}
\|P_{s}f\|_{B_{p,q}^{\alpha}}&=&\|P_{s}f\|_{p,\gamma_{d}}+\big(\int_{0}^{+\infty}\big(t^{k-\alpha}\|\frac{\partial^{k}P_{t}\big(P_{s}f\big)}{\partial t^{k}}\|_{p,\gamma_{d}}\big)^{q}\frac{dt}{t}\big)^{1/q}\\
&\leq&\|f\|_{p,\gamma_{d}}+\big(\int_{0}^{+\infty}\big(t^{k-\alpha}\|\frac{\partial^{k}P_{t}f}{\partial
t^{k}}\|_{p,\gamma_{d}}\big)^{q}\frac{dt}{t}\big)^{1/q}=\|f\|_{B_{p,q}^{\alpha}}.
\end{eqnarray*}

\item[ii)]  Now let us see that $\mathcal{J}_\beta $ is bounded on
$B_{p,q}^{\alpha}(\gamma_{d})$.
 Using the
Lebesgue's dominated convergence theorem and Minkowski's integral
inequality and  Jensen's inequality, we have
\begin{eqnarray*}
\|\frac{\partial^{k}P_{t}}{\partial t^{k}}\big(\mathcal{J}_\beta f
\big)\|_{p,\gamma_{d}}^{q}
&=&(\int_{\mathbb{R}^{d}}|\frac{\partial^{k} P_{t}}{\partial
t^{k}}(\frac{1}{\Gamma(\beta)}\int_{0}^{+\infty}s^{\beta} e^{-s}P_{s}f(x)\frac{ds}{s})|^{p}\gamma_{d}(dx))^{\frac{q}{p}}\\
&\leq&(\frac{1}{\Gamma(\beta)}\int_{0}^{+\infty}s^{\beta}e^{-s}(\int_{\mathbb{R}^{d}}|\frac{\partial^{k} P_{t}P_{s}f(x)}{\partial t^{k}}|^{p}\gamma_{d}(dx))^{\frac{1}{p}}\frac{ds}{s})^q\\
&\leq&\frac{1}{\Gamma(\beta)}\int_{0}^{+\infty}s^{\beta}e^{-s}\|\frac{\partial^{k}P_{t}P_{s}f}{\partial
t^{k}}\|_{p,\gamma_{d}}^{q}\frac{ds}{s},
\end{eqnarray*}
and then using Tonelli's Theorem,
\begin{eqnarray*}
\int_{0}^{+\infty}\big(t^{k-\alpha}\|\frac{\partial^{k}
P_{t}}{\partial t^{k}}\big(\mathcal{J}_\beta
f\big)\|_{p,\gamma_{d}}\big)^{q}\frac{dt}{t} & \leq&\frac{1}{\Gamma(\beta)}\int_{0}^{+\infty}s^{\beta}e^{-s}(\int_{0}^{+\infty}\big(t^{k-\alpha}\|\frac{\partial^{k} P_{t}\big(P_{s}f\big)}{\partial t^{k}}\|_{p,\gamma_{d}}\big)^{q}\frac{dt}{t})\frac{ds}{s}\\
&\leq&\frac{1}{\Gamma(\beta)}\int_{0}^{+\infty}s^{\beta}e^{-s}(\int_{0}^{+\infty}\big(t^{k-\alpha}\|\frac{\partial^{k}P_{t}f}{\partial t^{k}}\|_{p,\gamma_{d}}\big)^{q}\frac{dt}{t})\frac{ds}{s}\\
&=&\int_{0}^{+\infty}\big(t^{k-\alpha}\|\frac{\partial^{k}P_{t}f}{\partial
t^{k}}\|_{p,\gamma_{d}}\big)^{q}\frac{dt}{t}.
\end{eqnarray*}
Therefore
\begin{eqnarray*}
\| \mathcal{J}_\beta f\|_{B_{p,q}^{\alpha}}&=&\| \mathcal{J}_\beta
f\|_{p,\gamma_{d}}+\int_{0}^{+\infty}\big(t^{k-\alpha}\|\frac{\partial^{k}P_{t}}{\partial t^{k}}\big(\mathcal{J}_\beta f \big)\|_{p,\gamma_{d}}\big)^{q}\frac{dt}{t}\\
&\leq&\|
f\|_{p,\gamma_{d}}+\int_{0}^{+\infty}\big(t^{k-\alpha}\|\frac{\partial^{k}
P_{t}f}{\partial
t^{k}}\|_{p,\gamma_{d}}\big)^{q}\frac{dt}{t}=\|f\|_{B_{p,q}^{\alpha}}.
\end{eqnarray*}

\item[iii)]  Let us prove the $F_{p,q}^{\alpha}$-continuity of $P_t$ for
any $t >0$, the proof for $T_t$ is totally analogous. By the
Lebesgue's dominated convergence theorem and Minkowski's integral
inequality, we have
\begin{eqnarray*}
\left( \int_0^{\infty} (s^{k-\alpha}
\left|\frac{\partial^{k}P_{t}(P_{s}g)}{\partial s^{k}}(x) \right|)
^{q}\frac{ds}{s} \right) ^{1/q} &=&\left( \int_0^{\infty} (s^{k-\alpha}
\left|\int_{\mathbb{R}^{d}}p(t,x,y)\frac{\partial^{k}P_{s}g(y)}{\partial
s^{k}}dy\right|)^{q}\frac{ds}{s} \right) ^{1/q}\\
&\leq&\int_{\mathbb{R}^{d}}p(t,x,y)\left( \int_0^{\infty}
(s^{k-\alpha} \left|\frac{\partial^{k}P_{s}g(y)}{\partial
s^{k}} \right|) ^{q}\frac{ds}{s} \right) ^{1/q}dy\\
&=&  P_{t} \big(\left(\int_0^{\infty}(s^{k-\alpha}
\left|\frac{\partial^{k}P_{s} g}{\partial s^{k}} \right
|)^{q}\frac{ds}{s} \right)^{1/q}\big)(x).
\end{eqnarray*}
Therefore, by the $L^p$ continuity of $P_t$ we get
\begin{eqnarray*}
\|\left( \int_0^{\infty} (s^{k-\alpha}
\left|\frac{\partial^{k} P_{s}(P_{t}g)}{\partial s^{k}} \right|) ^{q}\frac{ds}{s} \right) ^{1/q}\|_{p,\gamma_{d}}&\leq&\|P_{t} \big(\left(\int_0^{\infty}(s^{k-\alpha} \left|\frac{\partial^{k} P_{s} g}{\partial s^{k}} \right |)^{q}\frac{ds}{s} \right)^{1/q}\big)\|_{p,\gamma_{d}}\\
&\leq&\|\left(\int_0^{\infty}(s^{k-\alpha}
\left|\frac{\partial^{k} P_{s} g}{\partial s^{k}} \right
|)^{q}\frac{ds}{s} \right)^{1/q}\|_{p,\gamma_{d}}
\end{eqnarray*}
Thus,
\begin{eqnarray*}
\|P_{t}g\|_{F_{p,q}^{\alpha}}&=&\|P_{t}g\|_{p,\gamma_{d}}+\|\left(
\int_0^{\infty}
(s^{k-\alpha} \left|\frac{\partial^{k} P_{s}(P_{t}g)}{\partial s^{k}} \right|) ^{q}\frac{ds}{s} \right) ^{1/q}\|_{p,\gamma_{d}}\\
&\leq&\|g\|_{p,\gamma_{d}}+\|\left(\int_0^{\infty}(s^{k-\alpha}
\left|\frac{\partial^{k} P_{s} g}{\partial s^{k}}
\right|)^{q}\frac{ds}{s} \right)^{1/q}\|_{p,\gamma_{d}}
=\|g\|_{F_{p,q}^{\alpha}}.
\end{eqnarray*}

\item[iv)] Now let us see that $\mathcal{J}_\beta $ is bounded on
$F_{p,q}^{\alpha}(\gamma_{d})$. By the Lebesgue's dominated
convergence theorem, Minkowski's integral inequality  and iii), we
have
\begin{eqnarray*}
&&\left(\int_0^{\infty}(s^{k-\alpha}
\left|\frac{\partial^{k}P_{s}}{\partial s^{k}}
\big(\mathcal{J}_{\beta}^{\gamma}g\big)(x) \right
|)^{q}\frac{ds}{s} \right)^{1/q}
\quad\quad \quad\quad \quad\\
&&\quad \quad \quad \quad \quad= \left( \int_0^{\infty} (s^{k-\alpha}
\left|\frac{\partial^{k}P_{s}}{\partial
s^{k}}\big(\frac{1}{\Gamma(\beta)}
\int_{0}^{+\infty}t^{\beta}e^{-t}P_{t}g(x)\frac{dt}{t}\big)\right|)^{q}\frac{ds}{s}\right)^{1/q}\\
&& \quad \quad \quad \quad \quad
\leq \frac{1}{\Gamma(\beta)}\int_{0}^{+\infty}t^{\beta}e^{-t}\left(
\int_0^{\infty} (s^{k-\alpha}
\left|\frac{\partial^{k}P_{s}(P_{t}g)}{\partial s^{k}}(x) \right|)
^{q}\frac{ds}{s} \right) ^{1/q}\frac{dt}{t},
\end{eqnarray*}
then, again by the Minkowski's integral inequality and  iii)
\begin{eqnarray*}
&&\|\left(\int_0^{\infty}(s^{k-\alpha}
\left|\frac{\partial^{k}P_{s}}{\partial s^{k}}
\big(\mathcal{J}_{\beta}^{\gamma}g\big) \right |)^{q}\frac{ds}{s}
\right)^{1/q}\|_{p,\gamma_{d}} \quad \quad\quad \quad \quad \quad
\quad\quad \quad\quad \quad\\
&&\quad \quad \quad \quad \quad \leq
\|\frac{1}{\Gamma(\beta)}\int_{0}^{+\infty}t^{\beta}e^{-t}\left(
\int_0^{\infty} (s^{k-\alpha} \left|\frac{\partial^{k}
P_{s}(P_{t}g)}{\partial s^{k}} \right|) ^{q}\frac{ds}{s} \right)
^{1/q}\frac{dt}{t}\|_{p,\gamma_{d}}\\
&&\quad \quad \quad \quad \quad \leq
\frac{1}{\Gamma(\beta)}\int_{0}^{+\infty} t^{\beta}e^{-t}\|\left(
\int_0^{\infty} (s^{k-\alpha} \left|\frac{\partial^{k}
P_{s}(P_{t}g)}{\partial s^{k}} \right|) ^{q}\frac{ds}{s} \right)
^{1/q}\|_{p,\gamma_{d}}\frac{dt}{t}\\
&&\quad \quad \quad \quad \quad
\leq\frac{1}{\Gamma(\beta)}\int_{0}^{+\infty}
t^{\beta}e^{-t}\|\left( \int_0^{\infty} (s^{k-\alpha}
\left|\frac{\partial^{k}P_{s} g}{\partial s^{k}} \right|)
^{q}\frac{ds}{s} \right)^{1/q}\|_{p,\gamma_{d}}\frac{dt}{t}\\
&&\quad \quad \quad \quad \quad  =\|\left( \int_0^{\infty}
(s^{k-\alpha} \left|\frac{\partial^{k} P_{s} g}{\partial s^{k}}
\right|) ^{q}\frac{ds}{s} \right) ^{1/q}\|_{p,\gamma_{d}}.
\end{eqnarray*}
Thus
\begin{eqnarray*}
\|\mathcal{J}_{\beta}^{\gamma}g\|_{F_{p,q}^{\alpha}}&=&\|\mathcal{J}_{\beta}^{\gamma}g\|_{p,\gamma_{d}}+
\|\left(\int_0^{\infty}(s^{k-\alpha}
\left|\frac{\partial^{k}P_{s}}{\partial s^{k}}
\big(\mathcal{J}_{\beta}^{\gamma}g\big) \right |)^{q}\frac{ds}{s}
\right)^{1/q}\|_{p,\gamma_{d}}
\\
&\leq&\|g\|_{p,\gamma_{d}}+\|\left( \int_0^{\infty} (s^{k-\alpha}
\left|\frac{\partial^{k}P_{s} g}{\partial s^{k}} \right|)
^{q}\frac{ds}{s} \right)
^{1/q}\|_{p,\gamma_{d}}=\|g\|_{F_{p,q}^{\alpha}}.
\end{eqnarray*}
 \end{enumerate}
\ep

Actually we can say more,
\begin{theorem}
Suppose that $\alpha\geq 0,\beta>0$. Then 
\begin{enumerate}
\item [i)] $\mathcal{J}_{\beta}$ is
bounded from $B_{p,q}^{\alpha}(\gamma_{d})$ to
$B_{p,q}^{\alpha+\beta}(\gamma_{d})$.
\item [ii)] $\mathcal{J}_{\beta}$ is
bounded from $F_{p,q}^{\alpha}(\gamma_{d})$ to
$F_{p,q}^{\alpha+\beta}(\gamma_{d})$.
\end{enumerate}
\end{theorem}
\dem 
\begin{enumerate}
\item [i)]  Let us denote $u(x,t)=P_{t}f(x)$ and
$U(x,t)=P_{t}\mathcal{J}_{\beta}f(x)$, using the representation of
$P_t$ (\ref{02}) we have,
$$U(x,t)=\displaystyle\int_{0}^{+\infty}T_{s}(\mathcal{J}_{\beta}f)(x)\mu_{t}^{(1/2)}(ds)$$
Therefore,
$$U(x,t_{1}+t_{2})=P_{t_{1}}(P_{t_{2}}(\mathcal{J}_{\beta}f))(x)=
\displaystyle\int_{0}^{+\infty}T_{s}(P_{t_{2}}(\mathcal{J}_{\beta}f))(x)\mu_{t_{1}}^{(1/2)}(ds).$$
Now, let $k,l$ be integer greater than $\alpha,\beta$
respectively, by differentiating $k$ times  respect to $t_{2}$ and $l$
times respect to $t_{1}$,

$$\displaystyle\frac{\partial^{k+l}U(x,t_{1}+t_{2})}{\partial
(t_{1}+t_{2})^{k+l}}=\int_{0}^{+\infty}T_{s}
(\frac{\partial^{k}P_{t_{2}}}{\partial
t_{2}^{k}}(\mathcal{J}_{\beta}f))(x)\frac{\partial^{l}}{\partial
t_{1}^{l}}\mu_{t_{1}}^{(1/2)}(ds).$$
Thus
$$\displaystyle\frac{\partial^{k+l} U(x,t)}{\partial
t^{k+l}}=\int_{0}^{+\infty}T_{s}(\frac{\partial^{k}P_{t_{2}}}{\partial
t_{2}^{k}}(\mathcal{J}_{\beta}f))(x)\frac{\partial^{l}}{\partial
t_{1}^{l}}\mu_{t_{1}}^{(1/2)}(ds),$$ if $t=t_{1}+t_{2}$ and
therefore, using the $L^p$ continuity of $T_s$ and (\ref{onesideineq})
\begin{eqnarray} \label{ineqU}
\nonumber \|\frac{\partial^{k+l} U(\cdot,t)}{\partial
t^{k+l}}\|_{p,\gamma}&\leq&
\int_{0}^{+\infty}\|T_{s} (\frac{\partial^{k}P_{t_{2}}}{\partial t_{2}^{k}}(\mathcal{J}_{\beta}f))\|_{p,\gamma} |\frac{\partial^{l}}{\partial t_{1}^{l}} \mu_{t_{1}}^{(1/2)}(ds)|\\
\nonumber &\leq&\int_{0}^{+\infty}\|
\frac{\partial^{k}P_{t_{2}}}{\partial
t_{2}^{k}}(\mathcal{J}_{\beta}f)\|_{p,\gamma}
|\frac{\partial^{l}}{\partial
t_{1}^{l}} \mu_{t_{1}}^{(1/2)}(ds)|\\
\nonumber &=&\|\frac{\partial^{k}P_{t_{2}}}{\partial t_{2}^{k}}(\mathcal{J}_{\beta}f)\|_{p,\gamma}\int_{0}^{+\infty}|\frac{\partial^{l}}{\partial t_{1}^{l}} \mu_{t_{1}}^{(1/2)}(ds)|\\
&\leq& C (t_{1})^{-l} \|\frac{\partial^{k}}{\partial
t_{2}^{k}}P_{t_{2}}\mathcal{J}_{\beta}f\|_{p,\gamma}
\end{eqnarray}

On  the other hand, using the representation of Bessel potential
(\ref{Beselrepre}) we have

$$P_{t}(\mathcal{J}_{\beta}f)(x)=\displaystyle\frac{1}{\Gamma(\beta)}\int_{0}^{+\infty}s^{\beta}e^{-s}P_{t+s}f(x)\frac{ds}{s}$$

then

\begin{eqnarray*}
\frac{\partial^{k}P_{t}}{\partial t^{k}}(\mathcal{J}_{\beta}f)(x)&=&\frac{1}{\Gamma(\beta)}\int_{0}^{+\infty}s^{\beta}e^{-s}\frac{\partial^{k} P_{t+s}f(x)}{\partial t^{k}}\frac{ds}{s}=\frac{1}{\Gamma(\beta)}\int_{0}^{+\infty}s^{\beta}e^{-s}\frac{\partial^{k}
P_{t+s}f(x)}{\partial (t+s)^{k}}\frac{ds}{s},
\end{eqnarray*}

and this implies that

\begin{eqnarray*}
\|\frac{\partial^{k} P_{t}}{\partial
t^{k}}(\mathcal{J}_{\beta}f)\|_{p,\gamma}&\leq&\frac{1}{\Gamma(\beta)}\int_{0}^{+\infty}s^{\beta}e^{-s}\|\frac{\partial^{k}P_{t+s}f}{\partial
(t+s)^{k}}\|_{p,\gamma}\frac{ds}{s},
\end{eqnarray*}
since $f\in B_{p,q}^{\alpha}(\gamma_{d})$. Now due to the fact
that the definition of $B_{p,q}^{\alpha}(\gamma_{d})$ is
independent on the integer $k>\alpha$ that we can choose, let us
take $k>\alpha+\beta$ and $l>\beta$, then
$k+l>\alpha+2\beta>\alpha+\beta$, this is  $k+l$ is an integer
greater than $\alpha+\beta$. Let us see now that

$$\displaystyle\big(\int_{0}^{+\infty}\big(t^{k+l-(\alpha+\beta)}\|\frac{\partial^{k+l} U(\cdot,t)}{\partial
t^{k+l}}\|_{p,\gamma}\big)^{q}\frac{dt}{t}\big)^{\frac{1}{q}}<+\infty.
$$ In fact, taking $t_1=t_2 = t/2$ in (\ref{ineqU}), we get
\begin{eqnarray*}
&&\displaystyle\big(\int_{0}^{+\infty}\big(t^{k+l-(\alpha+\beta)}\|\frac{\partial^{k+l}
U(\cdot,t)}{\partial
t^{k+l}}\|_{p,\gamma}\big)^{q}\frac{dt}{t}\big)^{\frac{1}{q}}
\quad \quad\quad \quad
\quad\\
&&\quad \quad\quad \quad \quad\leq C
\big(\int_{0}^{+\infty}\big(t^{k+l-(\alpha+\beta)}\|\frac{\partial^{k}P_{\frac{t}{2}}}{\partial
(\frac{t}{2})^{k}}(\mathcal{J}_{\beta}f)\|_{p,\gamma}(\frac{t}{2})^{-l}\big)^{q}\frac{dt}{t}\big)^{\frac{1}{q}}\\
&&\quad \quad\quad \quad \quad\leq
\frac{C}{\Gamma(\beta)}\big(\int_{0}^{+\infty}\big(t^{k-(\alpha+\beta)}\big(\int_{0}^{+\infty}s^{\beta}e^{-s}
\|\frac{\partial^{k}P_{s+\frac{t}{2}}f}{\partial(\frac{t}{2})^{k}}\|_{p,\gamma}\frac{ds}{s}\big)\big)^{q}\frac{dt}{t}\big)^{\frac{1}{q}}\\
&&\quad \quad\quad \quad
\quad\leq\frac{C}{\Gamma(\beta)}\big[\int_{0}^{+\infty}t^{(k-(\alpha+\beta))q}\big(\int_{0}^{t}s^{\beta}
\|\frac{\partial^{k}P_{s+\frac{t}{2}}f}{\partial(s+\frac{t}{2})^{k}}\|_{p,\gamma}\frac{ds}{s}\big)^{q}\\
&&\quad \quad\quad \quad \quad \quad \quad \quad\quad \quad \quad
\quad \quad + \big(\int_{t}^{+\infty}s^{\beta}
\|\frac{\partial^{k}P_{s+\frac{t}{2}}f}{\partial(s+\frac{t}{2})^{k}}\|_{p,\gamma}\frac{ds}{s}\big)^{q}\frac{dt}{t}\big]^{\frac{1}{q}}.
\end{eqnarray*}
Using again that $(a+b)^{q}\leq C_{q}(a^{q}+b^{q})$ if $a,b\geq
0,q\geq 1$, but since $(a+b)^{1/q}\leq a^{1/q}+b^{1/q}$ if $a,b\geq 0,q\geq
1$,
\begin{eqnarray*}
&&\frac{C}{\Gamma(\beta)}\big[\int_{0}^{+\infty}t^{(k-(\alpha+\beta))q}\big(\int_{0}^{t}s^{\beta}
\|\frac{\partial^{k}P_{s+\frac{t}{2}}f}{\partial(s+\frac{t}{2})^{k}}\|_{p,\gamma}\frac{ds}{s}\big)^{q}\quad
\quad\quad \quad
\quad\\
&&\quad \quad\quad \quad \quad \quad \quad +
\big(\int_{t}^{+\infty}s^{\beta}
\|\frac{\partial^{k}P_{s+\frac{t}{2}}f}{\partial(s+\frac{t}{2})^{k}}\|_{p,\gamma}\frac{ds}{s}\big)^{q}\frac{dt}{t}\big]^{\frac{1}{q}}\\
&&\quad \quad\quad \quad
\quad\leq\frac{C}{\Gamma(\beta)}\big[\int_{0}^{+\infty}t^{(k-(\alpha+\beta))q}\big(\int_{0}^{t}s^{\beta}
\|\frac{\partial^{k}P_{s+\frac{t}{2}}f}{\partial(s+\frac{t}{2})^{k}}\|_{p,\gamma}\frac{ds}{s}\big)^{q} \frac{dt}{t}\big]^{1/q}\\
&&\quad \quad \quad \quad \quad \quad \quad + \frac{C}{\Gamma(\beta)}\big[ \int_{0}^{+\infty}t^{(k-(\alpha+\beta))q}\big(\int_{t}^{+\infty}s^{\beta} \|\frac{\partial^{k}P_{s+\frac{t}{2}}f}{\partial(s+\frac{t}{2})^{k}}\|_{p,\gamma}\frac{ds}{s}\big)^{q}\frac{dt}{t}\big]^{\frac{1}{q}}\\
&&\quad \quad\quad \quad \quad=I+II.
\end{eqnarray*}

Now, using lemma \ref{kdecay} and since $\beta>0$
\begin{eqnarray*}
I&=&\frac{C}{\Gamma(\beta)}\big[\int_{0}^{+\infty}t^{(k-(\alpha+\beta))q}\big(\int_{0}^{t}s^{\beta}
\|\frac{\partial^{k}P_{s+\frac{t}{2}}f}{\partial(s+\frac{t}{2})^{k}}\|_{p,\gamma}\frac{ds}{s}\big)^{q}\frac{dt}{t}\big]^{\frac{1}{q}}\\
&\leq&\frac{C}{\Gamma(\beta)}\big[\int_{0}^{+\infty}t^{(k-(\alpha+\beta))q}\big(\int_{0}^{t}s^{\beta}
\|\frac{\partial^{k}P_{\frac{t}{2}}f}{\partial(\frac{t}{2})^{k}}\|_{p,\gamma}\frac{ds}{s}\big)^{q}\frac{dt}{t}\big]^{\frac{1}{q}}\\
&=&\frac{C}{\beta\Gamma(\beta)}\big(\int_{0}^{+\infty}\big(t^{k-\alpha}\|\frac{\partial^{k}P_{\frac{t}{2}}f}{\partial(\frac{t}{2})^{k}}\|_{p,\gamma}\big)^{q}\frac{dt}{t}\big)^{\frac{1}{q}}\\
&=&C_{\alpha,\beta}\big(\int_{0}^{+\infty}\big(u^{k-\alpha}\|\frac{\partial^{k}P_{u}f}{\partial
u^{k}}\|_{p,\gamma}\big)^{q}\frac{du}{u}\big)^{\frac{1}{q}}<+\infty,
\end{eqnarray*}
since  $f\in B_{p}^{\alpha,q}(\gamma_{d})$.

On the other hand, using the Hardy inequality, since
$k>\alpha+\beta$ and lemma \ref{kdecay} we get
\begin{eqnarray*}
II&=&\frac{C}{\Gamma(\beta)}\big(\int_{0}^{+\infty}t^{(k-(\alpha+\beta))q}\big(\int_{t}^{+\infty}s^{\beta}
\|\frac{\partial^{k}P_{s+\frac{t}{2}}f}{\partial(s+\frac{t}{2})^{k}}\|_{p,\gamma}\frac{ds}{s}\big)^{q}\frac{dt}{t}\big)^{\frac{1}{q}}\\
&\leq&\frac{C}{\Gamma(\beta)}\big(\int_{0}^{+\infty}t^{(k-(\alpha+\beta))q}\big(\int_{t}^{+\infty}s^{\beta}
\|\frac{\partial^{k}P_{s}f}{\partial s^{k}}\|_{p,\gamma}\frac{ds}{s}\big)^{q}\frac{dt}{t}\big)^{\frac{1}{q}}\\
&\leq&\frac{C}{\Gamma(\beta)}\frac{1}{k-(\alpha+\beta)}\int_{0}^{+\infty}\big(s^{k-\alpha}
\|\frac{\partial^{k}}{\partial
s^{k}}P_{s}f\|_{p,\gamma}\big)^{q}\frac{ds}{s}\big)^{\frac{1}{q}}<+\infty
\end{eqnarray*}

since $f\in B_{p,q}^{\alpha}(\gamma_{d})$. Thus
$\mathcal{J}_{\beta}f\in B_{p,q}^{\alpha+\beta}(\gamma_{d})$ and moreover,

$$\displaystyle\|\mathcal{J}_{\beta}f\|_{B_{p,q}^{\alpha+\beta}}\leq
C _{\alpha,\beta}\|f\|_{B_{p,q}^{\alpha}}.$$ 
\item [ii)] Let  $k>\alpha+\beta+1$  a fixed integer, let
 $f\in F_{p,q}^{\alpha}(\gamma_{d})$ and let us consider 
 $h={\mathcal{J}}_{\beta}f$.
 
 Let us consider two cases:
\begin{enumerate}
\item  If $\beta\geq 1$. By the change of variable  $u=t+s$ and using Hardy's inequality, we get

\begin{eqnarray*}
&&\big(\int_{0}^{+\infty}\big(t^{k-(\alpha+\beta)}|\frac{\partial^{k}P_{t}h(x)}{\partial
t^{k}}|\big)^{q}\frac{dt}{t}\big)^{1/q} \quad\quad \quad \quad \quad
\quad\quad \quad\quad \quad\\
&&\quad \quad \quad \quad \quad \leq \frac{1}{\Gamma(\beta)}\big(\int_{0}^{+\infty}t^{q(k-(\alpha+\beta))}\big(\int_{0}^{+\infty}s^{\beta}e^{-s}
|\frac{\partial^{k}P_{t+s}f(x)}{\partial(t+s)^{k}}|\frac{ds}{s}\big)^{q}\frac{dt}{t}\big)^{\frac{1}{q}}\\
&&\quad \quad \quad \quad \quad \leq \frac{1}{\Gamma(\beta)}\big(\int_{0}^{+\infty}t^{q(k-(\alpha+\beta))}\big(\int_{t}^{+\infty}(u-t)^{\beta-1}
|\frac{\partial^{k}P_{u}f(x)}{\partial u^{k}}|du\big)^{q}\frac{dt}{t}\big)^{\frac{1}{q}}\\
&&\quad \quad \quad \quad \quad \leq\frac{1}{\Gamma(\beta)}\big(\int_{0}^{+\infty}\big(\int_{t}^{+\infty}u^{\beta-1}
|\frac{\partial^{k}P_{u}f(x)}{\partial
u^{k}}|du\big)^{q}t^{q(k-(\alpha+\beta))-1}dt\big)^{\frac{1}{q}}\\
&&\quad \quad \quad \quad \quad \leq \frac{1}{\Gamma(\beta)}\frac{1}{k-(\alpha+\beta)}
\big(\int_{0}^{+\infty}\big(u^{k-\alpha}
|\frac{\partial^{k}P_{u}f(x)}{\partial
u^{k}}|\big)^{q}\frac{du}{u}\big)^{\frac{1}{q}}.
\end{eqnarray*}

Therefore,
\begin{eqnarray*}
&&\|\big(\int_{0}^{+\infty}\big(t^{k-(\alpha+\beta)}|\frac{\partial^{k}P_{t}h}{\partial
t^{k}}|\big)^{q}\frac{dt}{t}\big)^{\frac{1}{q}}\|_{p,\gamma}\quad \quad\quad \quad \quad \quad
\quad\quad \quad\quad \quad\\
&&\quad \quad \quad \quad \quad \leq
\frac{1}{\Gamma(\beta)(k-(\alpha+\beta))}
\|\big(\int_{0}^{+\infty}\big(u^{k-\alpha}
|\frac{\partial^{k}P_{u}f}{\partial
u^{k}}|\big)^{q}\frac{du}{u}\big)^{\frac{1}{q}}\|_{p,\gamma} <\infty,
\end{eqnarray*}
since $f\in  F_{p,q}^{\alpha}(\gamma_d)$. Thus ${\mathcal{J}}_{\beta}f\in  F_{p,q}^{\alpha+\beta}(\gamma_d)$.

\item  If $0<\beta< 1$.
 \begin{eqnarray*}
&&\big(\int_{0}^{+\infty}\big(t^{k-(\alpha+\beta)}|\frac{\partial^{k}P_{t} h(x)}{\partial
t^{k}}|\big)^{q}\frac{dt}{t}\big)^{\frac{1}{q}} \quad \quad\quad \quad \quad \quad
\quad\quad \quad\quad \quad\\
&&\quad \quad \quad \quad \quad \leq\frac{1}{\Gamma(\beta)}\big(\int_{0}^{+\infty}t^{q(k-(\alpha+\beta))}\big(\int_{0}^{+\infty}s^{\beta}e^{-s}
|\frac{\partial^{k}P_{t+s}f(x)}{\partial(t+s)^{k}}|\frac{ds}{s}\big)^{q}\frac{dt}{t}\big)^{\frac{1}{q}} \\
&&\quad \quad \quad \quad \quad \leq 
\frac{C}{\Gamma(\beta)}(\int_{0}^{+\infty}t^{q(k-(\alpha+\beta))-1}\big(\int_{0}^{t}s^{\alpha-1}e^{-s}
|\frac{\partial^{k}P_{t+s}f(x)}{\partial(t+s)^{k}}|ds\big)^{q}dt\big)^{\frac{1}{q}}\\
&&\quad \quad \quad \quad \quad+\frac{C}{\Gamma(\beta)}
(\int_{0}^{+\infty}t^{q(k-(\alpha+\beta))-1}\big(\int_{t}^{+\infty}s^{\alpha-1}e^{-s}
|\frac{\partial^{k}P_{t+s}f(x)}{\partial(t+s)^{k}}|ds\big)^{q}dt\big)^{\frac{1}{q}}\\
&&\quad \quad \quad \quad \quad=I + II.
\end{eqnarray*}
Now, $e^{-s}<1$ and as $\beta<1$, then $s^{\beta-1}<t^{\beta-1}$ for $t<s$. 

Hence again  by the change of variable  $u=t+s$ and using Hardy's inequality, we get
\begin{eqnarray*}
  II &\leq& \frac{C}{\Gamma(\beta)} \big(\int_{0}^{+\infty}t^{q(k-\alpha-1)-1}\big(\int_{t}^{+\infty}
|\frac{\partial^{k}P_{t+s}f(x)}{\partial(t+s)^{k}}|ds\big)^{q}dt\big)^{\frac{1}{q}}\\
&\leq&\frac{C}{\Gamma(\beta)} \big(\int_{0}^{+\infty}t^{q(k-\alpha-1)-1}\big(\int_{t}^{+\infty}
|\frac{\partial^{k}P_{u}f(x)}{\partial
u^{k}}|du\big)^{q}dt\big)^{\frac{1}{q}}\\
&\leq& \frac{C}{\Gamma(\beta)}\big(\int_{0}^{+\infty}\big(u^{k-\alpha}|\frac{\partial^{k}P_{u}f(x)}{\partial
u^{k}}|\big)^{q}\frac{du}{u}\big)^{\frac{1}{q}}
\end{eqnarray*}
On the other hand, using again  $e^{-s}<1$,
\begin{eqnarray*}
I^q &\leq&\frac{C}{\Gamma(\beta)}\int_{0}^{+\infty}t^{q(k-(\alpha+\beta))-1}\big(\int_{0}^{t}s^{\beta-1}
|\frac{\partial^{k}P_{t+s}f(x)}{\partial(t+s)^{k}}|ds\big)^{q}dt\\
&=&\frac{C}{\Gamma(\beta)\beta^{q}}\int_{0}^{+\infty}t^{q(k-\alpha)-1}\big(\frac{\beta
}{ t^{\beta}}\int_{0}^{t}s^{\beta-1}
|\frac{\partial^{k}P_{t+s}f(x)}{\partial(t+s)^{k}}|ds\big)^{q}dt
\end{eqnarray*}
Now, as $\beta >0$,
$\displaystyle\int_{0}^{t}s^{\beta-1}ds=\frac{t^{\beta}}{\beta},$ then
using Jensen's inequality for the probability measure 
$\displaystyle\frac{\beta}{t^{\beta}}s^{\beta-1} ds$ and Fubini's theorem
\begin{eqnarray*}
I^q&\leq& \frac{C}{\Gamma(\beta)\beta^{q}}\int_{0}^{+\infty}t^{q(k-\alpha)-1}\big(\frac{\beta
}{ t^{\beta}}\int_{0}^{t}s^{\beta-1}
|\frac{\partial^{k}P_{t+s}f(x)}{\partial(t+s)^{k}}|^{q}ds\big)dt\\
&=& \frac{C}{\Gamma(\beta)\beta^{q-1}}\int_{0}^{+\infty}s^{\beta-1}\big(\int_{s}^{+\infty}t^{q(k-\alpha)-\beta-1}
|\frac{\partial^{k}P_{t+s}f(x)}{\partial(t+s)^{k}}|^{q}dt\big)ds\\
&\leq& \frac{C}{\Gamma(\beta)\beta^{q-1}}\int_{0}^{+\infty}s^{\beta-1}\big(\int_{s}^{+\infty}(t+s)^{q(k-\alpha)-\beta-1}
|\frac{\partial^{k}P_{t+s}f(x)}{\partial(t+s)^{k}}|^{q}dt\big)ds
\end{eqnarray*}
as $q(k-\alpha)-\beta-1>0$, since $0<\beta <1$. Finally, again by  the change of variable  $u=t+s$ and using Hardy's inequality, we get
\begin{eqnarray*}
I^q&\leq&\frac{C}{\Gamma(\beta)\beta^{q-1}}\int_{0}^{+\infty}s^{\beta-1}\big(\int_{2s}^{+\infty}u^{q(k-\alpha)-\beta-1}
|\frac{\partial^{k}P_{u}f(x)}{\partial u^{k}}|^{q}du\big)ds\\
&\leq&\frac{C}{\Gamma(\beta)\beta^{q-1}}\int_{0}^{+\infty}s^{\beta-1}\big(\int_{s}^{+\infty}u^{q(k-\alpha)-\beta-1}
|\frac{\partial^{k}P_{u}f(x)}{\partial u^{k}}|^{q}du\big)ds\\
&\leq&\frac{C}{\Gamma(\beta)\beta^{q-1}}\int_{0}^{+\infty}\big(u^{k-\alpha}|\frac{\partial^{k}P_{u}f(x)}{\partial
u^{k}}|\big)^{q}\frac{du}{u}.
\end{eqnarray*}
Therefore
\begin{eqnarray*}
\|\big(\int_{0}^{+\infty}\big(t^{k-(\alpha+\beta)}|\frac{\partial^{k}P_{t}h}{\partial
t^{k}}|\big)^{q}\frac{dt}{t}\big)^{\frac{1}{q}}\|_{p,\gamma} \leq
C_{k,\alpha,\beta }\|\big(\int_{0}^{+\infty}\big(u^{k-\alpha}|\frac{\partial^{k}P_{u}f}{\partial
u^{k}}|\big)^{q}\frac{du}{u} \big)^{\frac{1}{q}}\|_{p,\gamma} < \infty.
\end{eqnarray*}

Thus ${\mathcal{J}}_{\beta}f\in  F_{p,q}^{\alpha+\beta}(\gamma_d)$, for $0< \beta<1$.

In both cases we have,
\begin{eqnarray*}
\|{\mathcal{J}}_{\beta}f\|_{F_{p,q}^{\alpha+\beta}}&=&\|{\mathcal{J}}_{\beta}f\|_{p,\gamma}+
\|\big(\int_{0}^{+\infty}\big(t^{k-(\alpha+\beta)}|\frac{\partial^{k}P_{t}{\mathcal{J}}_{\beta}f}{\partial
t^{k}}|\big)^{q}\frac{dt}{t}\big)^{\frac{1}{q}}\|_{p,\gamma} \\
&\leq &C_{\beta} \|f\|_{p,\gamma}+C_{k,\alpha,\beta}\|\big(\int_{0}^{+\infty}\big(u^{k-\alpha}|\frac{\partial^{k}P_{u}}{\partial u^{k}}|\big)^{q}\frac{du}{u} \big)^{\frac{1}{q}}\|_{p,\gamma} \\
&\leq& C_{k,\alpha,\beta }\|f\|_{F_{p,q}^{\alpha}}.
\end{eqnarray*}  \ep
\end{enumerate}
\end{enumerate}


\begin{thebibliography}{99}

\bibitem{berg} Berg C., Reus C. J.P. {\it Density questions in the classical theory of moments.}
Annales de l'institut Fourier, tome 31, 3 (1981), 99-114.
\bibitem{ButzBer}
 Butzer, P. L. \& Berens, H. \emph{Semi-groups of operators and approximation.} Die Grundkehren der mathematischen Wissenschaften, Ban145. Springer-Verlag. New York, 1967.
\bibitem{ep}
  Epperson, J. \emph{Triebel-Lizorkin spaces for Hermite expansions.} Studia Math, {\bf 114} (1995),199--209.
\bibitem{fm91}
  Frazier M, Jawerth B, Weiss G. \emph{Littlewood Paley Theory and the Study of Functions Spaces.} CBMS-Conference Lecture Notes 79.
 Amer. Math. Soc. Providence RI , 1991.
 \bibitem{forscotur}
Forzani, L., Scotto, R, and Urbina, W{.} \emph{Riesz and Bessel
Potentials, the $g_k$ functions and an Area function, for  the
Gaussian measure $\gamma_d$.} Revista de la Uni\'on Matem\'atica
 Argentina (UMA), vol 42 (2000), no.1,17--37.
 \bibitem{forur}
  Forzani, L. \& Urbina,W. {\em Poisson-Hermite representation of solutions of the equation
$ \frac{\partial^2}{\partial t^2}u(x,t) + \Delta_x u(x,t) - 2x\cdot \nabla_x u(x,t) = 0$. } Proceedings 5th International Conference on
Approximation and Optimization in the Caribbean. Approximation, Optimization and Mathematical 
Economics (2001). 109-115. Springer Verlag.
\bibitem{gs96}
 Gatto A. E, Segovia C, V\'{a}gi S. \emph{On Fractional Differentiation and Integration on Spaces of Homogeneous Type.} Rev. Mat. Iberoamericana,  {\bf 12} (1996), 111--145.
 \bibitem{gs99}
 Gatto A. E, V\'{a}gi S. \emph{On Sobolev Spaces of Fractional Order and $\epsilon-$families of Operators on Spaces of Homogeneous Type.} Studia Math. {\bf 133} (1999), no.1, 19--27.
\bibitem{gu}
Guti\'errez, C.  {\em On the Riesz transforms 
for the Gaussian measure.} J. Func. Anal. 120 (1) (1994) 107-134
\bibitem{lu}
L\'{o}pez I{.} and  Urbina, W{.} \emph{Fractional Differentiation
for the Gaussian Measure and Applications. } Bull. Sciences Math, 2004,
{\bf 128}, 587--603.
\bibitem{iris}
L\'{o}pez, I{.} \emph{Introduction to the Besov Spaces and
Triebel-Lizorkin Spaces for Hermite and Laguerre expansions and
some applications.} J. of Math \& Stat. {\bf 1} (2005), no. 3, 172-179.
 \bibitem{ma}
Mazet, O. {\em Classification des demi-groupes de diffusion sur
$\mathbb{R}$ asocie a une famille de polynomes orthogonaux}.
S\'eminaire de Probabilit\'es, Lecture Notes in Math, vol.
1655, 40-54. Springer-Verlag (1997).
\bibitem{me}
 Meyer, P.A. \emph{Transformations de Riesz pour le lois Gaussiens.}  Lecture Notes in Math 1059   1984, Springer-Verlag 179--193.
 \bibitem{ebner} Pineda, E. \emph{Algunos T\'opicos en An\'alisis Arm\'onico Gaussiano: Comportamiento en la frontera,  Espacios de Besov-Lipstchitz y de Triebel-Lizorkin}.
 Doctoral Thesis, Facultad de Ciencias, UCV, Caracas. In preparation.
\bibitem{sp97}
  Sj\"{o}gren P. \emph{Operators associated with the Hermite semigroup- a survey.} J. Fourier Anal. Appl. {\bf 3} (1997), Special Issue, 813--823.
\bibitem{se70}
  Stein E. \emph{Singular integrals and differentiability properties of functions} Princeton Univ. Press. Princeton, New Jersey, 1970.
   \bibitem{trie}
 Triebel, H  \emph{Interpolation theory, function spaces differential operators.} Noth Holland, 1978.
 \bibitem{trie1}
 Triebel, H  \emph{Theory of function spaces.} Birkh\"auser Verlag, Basel, 1983.
   \bibitem{trie2}
 Triebel, H  \emph{Theory of function spaces II.} Birkh\"auser Verlag, Basel, 1992.
\bibitem{uw98}
  Urbina W. \emph{An\'{a}lisis Arm\'{o}nico Gaussiano: una visi\'{o}n panor\'{a}mica}. Trabajo de Ascenso, Facultad de Ciencias, UCV, Caracas, 1998.
  \bibitem{wat}
 Watanabe, S  \emph{Stochastic Differential Equations and Malliavin Calculus}, Tata Institute of Fundamental Research,Springer Verlag, Berlin,1984.
\end{thebibliography}
\end{document}